\documentclass[10pt,preprint]{elsarticle}
\usepackage[english]{babel}
\usepackage[utf8]{inputenc}
\usepackage[T1]{fontenc}
\usepackage[acronyms]{glossaries}
\usepackage{writemaths}
\usepackage{paralist}
\usepackage{enumitem}
\usepackage{subfig}
\usepackage{xcolor}
\usepackage{ifthen}

\usepackage{fullpage}

\makeglossaries

\newacronym{coc}{CoC}{\textit{change of coordinates}}
\newacronym{com}{CoM}{\textit{change of measure}}
\newacronym{st}{ST}{\textit{seismic tomography}}
\newacronym{td}{TD}{\textit{transient diffusion}}
\newacronym{kl}{KL}{\textit{Karhunen--Loève}}
\newacronym{se}{SE}{squared exponential}
\newacronym{mc}{MC}{\textit{Monte--Carlo}}
\newacronym{mcmc}{MCMC}{\textit{Markov Chain Monte--Carlo}}
\newacronym{mh}{MH}{\textit{Metropolis--Hastings}}
\newacronym{rrmse}{RRMSE}{\textit{relative root mean squared error}}
\newacronym{rmsre}{RMSRE}{\textit{root mean squared relative error}}
\newacronym{svd}{SVD}{\textit{singular value decomposition}}
\newacronym{psp}{PSP}{\textit{pseudo spectral projection}}
\newacronym{map}{MAP}{\textit{maximum~{a posteriori}}}
\newacronym{iid}{i.i.d.}{\textit{independent and identically distributed}}
\newacronym{sg}{SG}{\textit{sparse grid}}
\newacronym{pc}{PC}{\textit{polynomial chaos}}
\newacronym{eof}{EOF}{\textit{empirical orthogonal function}}
\newacronym{ess}{ESS}{\textit{effective sample size}}

\journal{Journal of Computational Physics}

\date{\today}

\begin{document}
\title{Change of Measure for Bayesian Field Inversion with Hierarchical Hyperparameters Sampling}

\author[cea,mines]{Nad\`ege Polette\corref{cor1}}
\ead{nadege.polette@minesparis.psl.eu}

\author[cnrs]{Olivier Le Ma\^itre}
\ead{olivier.le-maitre@polytechnique.edu}

\author[cea]{Pierre Sochala}
\ead{pierre.sochala@cea.fr}

\author[mines]{Alexandrine Gesret}
\ead{alexandrine.gesret@minesparis.psl.eu}

\cortext[cor1]{Corresponding author}
\address[cea]{CEA, DAM, DIF, F-91297 Arpajon, France}
\address[mines]{Mines Paris PSL, Geosciences center, 77300 Fontainebleau, France}
\address[cnrs]{CMAP, CNRS, Inria, Ecole Polytechnique, IPP, 91120 Palaiseau, France}

\begin{frontmatter}

\begin{abstract} 
\noindent   
This paper proposes an effective treatment of hyperparameters in the Bayesian inference of a scalar field from indirect observations. 
Obtaining the joint posterior distribution of the field and its hyperparameters is challenging. The infinite dimensionality of the field requires a finite parametrization that usually involves hyperparameters to reflect the limited prior knowledge. In the present work, we consider a Karhunen-Loève (KL) decomposition for the random field and hyperparameters to account for the lack of prior knowledge of its autocovariance function. The hyperparameters must be inferred. To efficiently sample jointly the KL coordinates of the field and the autocovariance hyperparameters, we introduce a change of measure to reformulate the joint posterior distribution into a hierarchical Bayesian form. The likelihood depends only on the field's coordinates in a fixed KL basis, with a prior conditioned on the hyperparameters. We exploit this structure to derive an efficient Markov Chain Monte Carlo (MCMC) sampling scheme based on an adapted Metropolis--Hasting algorithm. We rely on surrogate models (Polynomial Chaos expansions) of the forward model predictions to further accelerate the MCMC sampling. A first application to a transient diffusion problem shows that our method is consistent with other approaches based on a change of coordinates (Sraj et al., 2016). A second application to a seismic traveltime tomography highlights the importance of inferring the hyperparameters. A third application to a $2$D anisotropic groundwater flow problem illustrates the method on a more complex geometry.
\end{abstract}


\end{frontmatter}

\section{Introduction}\label{sec:intro}
Inverse problems arise in many situations whenever one searches for information about a physical system based on observations~\cite{tarantola2005}. The Bayesian approach~\cite{sivia2006} is widely used to solve inverse problems in a probabilistic framework and has been applied in various geophysical applications~\cite{malinverno2002,sochala2021a,grana2022,rounce2020,ying2019}. The Bayesian approach is attractive because it requires a weak~\textit{a priori} on the unknown parameters and provides a full estimation of the parameter distributions. In practice, \gls{mcmc} methods provide algorithms to sample from the posterior distribution~\cite{doucet2013}. \gls{mcmc} methods require a considerable number of (generally expensive) forward model evaluations, and the convergence of the sampling might be difficult to reach, especially when dealing with a high-dimensional parameter space. The computational cost induced by the sampling can be reduced using surrogate models, where the forward model predictions are replaced by fast to evaluate approximations~\cite{marzouk2007a,navarro2018,sochala2021a}.

In this study, we are interested in the estimation of a scalar field through a set of indirect and noisy observations. Such inference problem is challenging because of the infinite dimensionality of the field. Previous works~\cite{marzouk2007a,marzouk2009} have already addressed this issue by means of \gls{kl} decomposition~\cite{karhunen1946,loeve1977,meyer2003} to obtain a finite dimensional parametrization of the field. The \gls{kl} expansion represents a random field using the eigenelements of its autocovariance function as a decomposition basis. The inference then consists in identifying the \gls{kl} coordinates in the basis of the dominant eigenmodes. In practice, the autocovariance function depends on hyperparameters $\bq$ that can be determined by optimizing a criterion such as the likelihood estimation~\cite{rasmussen2005}, the leave-one-out error or by parametric estimation~\cite{uribe2020,meles2022}. Nevertheless, fixing a value of the hyperparameters can lead to overconfident results. This motivated the development of methods that jointly explore the parameters and hyperparameters spaces during the inverse problem solving. For stationnary fields, it is straightforward to deal with variable variances, as it reduces to a scaling of the \gls{kl} coordinates~\cite{marzouk2009,khatoon2023}. On the contrary, the other hyperparameters, like correlation length and exponent, are more delicate to infer since they affect directly the \gls{kl} basis, raising significant computational challenges~\cite{laloy2013}. Several methods~\cite{tagade2014,sraj2016,latz2019} have been proposed to deal with the computational cost induced by the eigenmodes basis $\bq$-dependency. In particular, \cite{sraj2016} introduced a method based on a fixed reference basis over the hyperparameters space. The reference basis is made of the dominant modes of the $\bq$-averaged autocovariance function. The dependency on the hyperparameters is transferred to the coordinates through a linear transformation called~\gls{coc}. This construction minimizes the $\bq$-averaged field representation error. Despite the benefit of the reference basis, the \gls{coc} method is limited to cases where the change of coordinates matrix is $\bq$-continuous~\cite{sraj2016,siripatana2020}. Indeed, the computation of this matrix is computationally expensive since it requires solving an eigenvalue problem at each \gls{mcmc} step. In addition, multiplicities in the eigenvalues and crossing of the eigenbranches with changes of the hyperparameters makes it difficult to ensure a smooth \gls{coc} as required by the \gls{mcmc} samplers.

The aim of this paper is to introduce a new method to alleviate the difficulties of the \gls{coc} method by transferring the $\bq$-dependency to the prior distribution of the field coordinates. In this approach, called \gls{com}, the coordinates of the field follow a Gaussian distribution whose covariance matrix depends smoothly on the hyperparameters and does not require to decompose the $\bq$-dependent autocovariance function at each \gls{mcmc} step. Another advantage of our formulation is that the posterior distribution only depends on $\bq$ through the prior distribution. To accelerate the \gls{mcmc} sampling, we rely on \gls{pc} surrogate models~\cite{wiener1938a,ghanem1991a,xiu2002} for both model predictions and \gls{com} derived quantities. 

The paper is structured as follows. Section~\ref{sec:bayes-inference} presents the Bayesian framework and the field parametrization. Section~\ref{sec:com} focuses on the derivation of the coordinates distribution in the \gls{com} method, especially its covariance matrix expression and the associated sampling procedure. Section~\ref{sec:acceleration} describes the \gls{pc} surrogate model used to replace the forward model, as well as the \gls{com} derived quantities. The \gls{com} method is compared with the \gls{coc} method on a transient diffusion case in Section~\ref{sec:transdiff} and applied to a traveltime seismic tomography problem in Section~\ref{sec:tomo}. Section~\ref{sec:groundwater} deals with the inference of a $2$D anisotropic random field in a groundwater flow application. Conclusions and perspectives are drawn in Section~\ref{sec:ccl}.

\section{Bayesian field inference}\label{sec:bayes-inference}
We detail the Bayesian framework in Section~\ref{sec:bayes-form}, and then we present the dimension reduction technique with the \gls{kl} decomposition depending on hyperparameters in Section~\ref{sec:representation-field}.

\subsection{Bayes' formulation}\label{sec:bayes-form}
We are interested in the inference of a scalar field $g \in \mathbb{G} \subset L^2(\Omega)$ defined over a compact spatial domain $\Omega\subset\R^d$, where $d \in \{1,2,3\}$ is the dimension of the physical space. 
Denoting $(\Theta, \mathfrak{S}, P)$ an abstract probability space and $\theta \in \Theta$ a particular random event, the field $g$ is assumed to be a particular realization of a random process $G$, that is $g(\bx) = G(\bx,\theta)$.

The definition of the inference problem relies on a set of $N$ indirect observations $\bs{d}_\mathrm{obs} \in \R^{N}$, an observation noise model, as well as a forward model $M: g \mapsto \bs{d}$ that predict $\bs{d}$ for a given field $g$, where $M$ is for instance a system of PDEs. This relation is not invertible in general and cannot be used to obtain a closed-form solution for $g$, such as $g\coloneqq M^{-1}\bs{d}_\mathrm{obs}$ composed with an observation operator. The inverse problem can be solved in several ways~\cite{tarantola2005, hansen2010} and we choose here the Bayesian approach, which characterizes the full posterior distribution of $g$.
According to the Bayes' formula, the posterior distribution of the field $g$ reads
\begin{equation}\label{eq:general-posterior}
\pi_\mathrm{post}(g| \bs{d}_\mathrm{obs}) \propto \mc{L}(\bs{d}_\mathrm{obs}|g)\pi_\mathbb{G}(g),
\end{equation}
where $\mc{L}(\bs{d}_\mathrm{obs}|g)$ denotes the likelihood of the observations given the field $g$ and $\pi_\mathbb{G}(g)$ its prior distribution. The posterior can be sampled with \gls{mcmc} methods. The two main difficulties to use the posterior distribution~\eqref{eq:general-posterior} in practice are the infinite dimension of the field and the likelihood computation cost, based on forward model evaluations at each \gls{mcmc} step. This motivates the use of a dimension reduction technique to approximate the field in a low dimensional space. Besides, the computational cost can be decreased with surrogate models that replace the forward model predictions.

\subsection{Dimension reduction of random fields}\label{sec:representation-field}
This section provides a parametrization of the field $g$ along with its associated sample space $\mathbb{G}$ and its prior probability distribution $\pi_\mathbb{G}$ for the Bayesian inference. Several approaches are possible to obtain finite dimensional representation of the field. The nodal representation relies on interpolation between values defined at a finite set of points constituting a grid. This representation is particularly expensive when a fine mesh is required to accomodate small scale features. Transdimensional approaches~\cite{bodin2012,pianaagostinetti2015,belhadj2018} have been developed to mitigate this issue by iteratively optimizing the number of mesh cells. These approches are however incompatible with the use of surrogate models since the field parametrization changes along the \gls{mcmc} steps. We choose instead a modal representation of the field~\cite{meyer2003,musolas2021,akian2022,zhang1994} which offers a low-dimensional parametrization that is fixed during the inference. 

Consider a random field $G(\bx,\theta)$ with mean $\mu:\Omega\mapsto\R$ and autocovariance function $k: \Omega^2 \mapsto \R$, defined as
\begin{equation}\label{eq:kernel-def}
\mu(\bx) \coloneqq \E_\Theta\left(G(\bx, \cdot)\right) \quad \text{and} \quad
k(\bx,\by) \coloneqq \E_\Theta\left( \left(G(\bx, \cdot)-\mu(\bx)\right) \left(G(\by,\cdot) -\mu(\by)\right) \right),
\end{equation}
where $\E_\Theta$ denotes the expectation. Without loss of generality, we set $\mu=0$ to alleviate notations.
The truncated \gls{kl} expansion $G^r$ of $G$ writes
\begin{equation}\label{eq:KL}
  G^r(\bx,\theta) \coloneqq \suml_{i=1}^{r} \lambda_i^{1/2}u_i(\bx)\eta_i(\theta), \quad\text{such that}\quad G(\bx,\theta) = \lim_{r\rightarrow \infty} G^r(\bx,\theta).
\end{equation}
The couples $\left(u_i, \lambda_i\right)_{i \in \N^{\ast}}$ are the eigenfunctions and associated eigenvalues (sorted in descending order) of the autocovariance function $k$, obtained by solving the Fredholm equation of the second kind,
\begin{equation}\label{eq:fredholm}
  \forall \bx \in \Omega, \ \forall i \in \N^{\ast}, \qquad \Int_\Omega k(\bx,\by)u_i(\by)d\by = \lambda_iu_i(\bx).
\end{equation}
The \gls{kl} decomposition is bi-orthonormal in the sense that the coordinates $\eta_i = \lambda^{-1/2}\scal{G(\theta)}{u_i}$ are uncorrelated with unitary variance and the eigenfunctions $u_i$ are orthonormal with respect to the inner product of $L^2(\Omega)$, 
\begin{equation}\label{eq:KL-orthonormality}
\E(\eta_i\eta_j)=\delta_{i,j}\quad\text{and}\quad\scal{u_i}{u_j}_\Omega \coloneqq \Int_\Omega u_i(\bx)u_j(\bx) d\bx = \delta_{i,j}.
\end{equation}
In addition, the truncated decomposition minimizes the representation error in the $L^2$-sense~\cite{meyer2003}: $\E_\theta \left(\norm{G-G^r}_{L^2(\Omega)}\right) = \suml_{i>r} \lambda_i$. The inference problem~\eqref{eq:general-posterior}  for $G^r$ can be recast into an identification of the $r$-dimensional coordinates vector $\bs{\eta}$,
\begin{equation}\label{eq:kl-posterior}
  \pi_\mathrm{post}(\bs{\eta}| \bs{d}_\mathrm{obs}) \propto \mc{L}(\bs{d}_\mathrm{obs}|\bs{\eta})\pi(\bs{\eta}).
\end{equation}
In the following, a Gaussian random process for the field is considered, and therefore $\bs{\eta}\sim\mathcal{N}(0,\mathrm{I}_r)$. 

The autocovariance function $k$ may depend on some hyperparameters $\bq \in \mathbb{H}$ that are poorly known~\cite{tagade2014,latz2019}. In case of an uncertain correlation length, one can choose a small value of it, use a \gls{kl} truncature with many modes and then apply model reduction techniques~\cite{cui2014,li2014}. Another solution consists in integrating the hyperparameters in the inference problem. This yields
\begin{equation}
  \pi_\mathrm{post}\left( \bs{\eta}, \bq \left| \bs{d}_\mathrm{obs}\right. \right) \propto \mc{L}(\bs{d}_\mathrm{obs} | \bs{\eta}, \bq)\pi(\bs{\eta}, \bq),
\end{equation}
where $\pi(\bs{\eta}, \bq)$ is the joint prior distribution of the \gls{kl} coordinates and hyperparameters.
The \gls{kl} expansion basis and the coordinates of the field $G$ depend on the hyperparameters,
\begin{equation}\label{eq:kl-qdep}
  G(\bx,\theta) \simeq \suml_{i=1}^{r} \ls{i}{\bq}u_i(\bx,\bq)\eta_i(\theta), \quad\text{with}\quad
  \eta_i(\theta) = \scal{\ls[-]{i}{\bq}u_i(\cdot,\bq)}{G(\cdot,\theta)}_\Omega.
\end{equation}
The estimation of the coordinates and hyperparameters in decomposition~\eqref{eq:kl-qdep} faces several difficulties, at the sampling stage as well as to construct surrogate models. First, the eigenelements must be computed at each \gls{mcmc} step. Further, the $\bq$-continuity of the field is not guaranteed due to  the non uniqueness of the eigenelements, particularly in the case of crossing eigenbranches~\cite{siripatana2020}. To get rid of this dependency in the field representation, \cite{sraj2016} introduces a fixed reference basis to expand the field. The method proposed in this work follows this idea by expressing the hyperparameter dependencies through the prior of the coordinates in the reference basis.

\section{Change of measure method}\label{sec:com}
Sections~\ref{subsec:ref-basis} and \ref{subsec:com-coc} briefly describe the reference basis and the change of coordinates method used to estimate the KL coordinates.
Section \ref{subsec:com} presents the novel approach, called \gls{com}, that overcomes the change of coordinates limitations. Section~\ref{subsec:com-implementation} details the log-posterior expression and the sampling strategy for the \gls{com}.

\subsection{Reference basis}\label{subsec:ref-basis}
The motivation of using a reference basis is to eliminate the eigenfunctions $\bq$-dependency. Such basis could be computed from the kernel $k(\cdot,\cdot,\ol{\bq})$ associated with a particular value of $\bq$, its mean value $\ol{\bq} = \E_\mathbb{H}(\bq)$ or from the kernel averaged over the hyperparameters domain $\mathbb{H}$. We opt for the latter choice because it minimizes the representation error of the field in average over $\bq$~\cite{sraj2016}. We denote $\ol{k}$ the averaged autocovariance function,
  \begin{equation}\label{eq:kref-def}
  \forall \bx, \by \in \Omega, \qquad \ol{k}(\bx,\by) \coloneqq \E_\mathbb{H}(k(\bx,\by,\cdot)) \coloneqq \Int_{\mathbb{H}} k(\bx, \by, \bq)\pi_\mathbb{H}(\bq)d\bq,
  \end{equation}
  where $\pi_\mathbb{H}$ denotes the hyperparameters prior.
  The reference eigenelements $\{\ol{u}_i, \ol{\lambda}_i\}_{i \in \N^{\ast}}$ are obtained by solving the eigenvalue problem for the reference kernel,
 \begin{equation}\label{eq:Bref-def}
 \forall \bx \in \Omega,\ \forall i \in \N^{\ast}, \qquad \Int_\Omega \ol{k}(\bx,\by)\ol{u}_i(\by)d\by = \ol{\lambda}_i\ol{u}_i(\bx).
 \end{equation}
The field approximation in the reference basis is 
\begin{equation}\label{eq:ref-decomposition}
 G(\bx,\theta) \simeq \ol{G}^r(\bx,\theta) \coloneqq \sum_{i=1}^r \ls{i}{r}\ol{u}_i(\bx)\xi_i(\theta),\quad\text{with}\quad\xi_i(\theta)=\scal{\ls[-]{i}{r}\ol{u}_i}{G(\cdot,\theta)}_\Omega
\end{equation}
The following sections propose two methods to use this parametrization of $G$ in the inference of the coordinates and hyperparameters joint distribution.

\subsection{Change of coordinates}\label{subsec:com-coc}
In~\cite{sraj2016}, the coordinates $\bs{\xi}$ are obtained by a transformation of the coordinates in the $\bq$-dependent basis defined previously in Eq.~\eqref{eq:kl-qdep},
\begin{equation}
G^r(\bx,\theta) = \suml_{i=1}^{r} \ls{i}{\bq}u_i(\bx,\bq)\eta_i(\theta) \simeq \suml_{i=1}^{r} \ls{i}{r} \ol{u}_i(\bx)\xi_i(\theta,\bq).
\end{equation}
The two sets of coordinates $\bxi(\theta,\bq)$ and $\bs{\eta}(\theta)$ are related by the change of coordinates matrix $B(\bq)$ defined by
\begin{equation}
 B(\bq)_{ij} = \ls[-]{i}{r} \scal{\ls{j}{\bq} u_j(\cdot,\bq)}{\ol{u}_i}_\Omega, \quad \bxi(\theta,\bq) = B(\bq)\bs{\eta}(\theta).
\end{equation}
The inference problem here writes 
\begin{equation}
 \pi_\mathrm{post}\left( \bs{\eta}, \bq \left| \bs{d}_\mathrm{obs}\right. \right) \propto \mc{L}(\bs{d}_\mathrm{obs} | \bs{\eta}, \bq)\pi(\bs{\eta})\pi_\mathbb{H}(\bq), \text{ with } \bs{\eta} \sim \mc{N}(0,\mathrm{I}_r).
\end{equation}
In this formulation, for fixed $\bs{\eta}$, changing the hyperparameters modifies the coordinates in the reference basis and thus the shape of the field. Therefore, changing the hyperparameter modifies the realization of the process, which is not so natural. In addition, the posterior probability distribution of the coordinates $\bs{\eta}$ depends on the hyperparameters: for each value of $\bq$, a different posterior distribution for $\bs{\eta}$ is sampled, which can complicate the inference process. 
Moreover, the dependency of $B$ on the eigenelements $\{u_i(\cdot,\bq),\lambda_i(\bq)\}_{1 \leq i \leq r}$ requires to solve the eigenvalue problem~\eqref{eq:fredholm} for each new value of $\bq$. In~\cite{sraj2016,siripatana2020} a surrogate model is built for $B(\bq)$. However, this strategy suffers from the numerical ambiguity of $B(\bq)$ due to the choice of the eigenvectors orientation and indexation~\cite{siripatana2020}. 
These drawbacks lead to reformulate the field parametrization thanks to a \acrfull{com}, where the hyperparameters values affects the prior of the proposed field and not the field values themselves. 

\subsection{Change of measure}\label{subsec:com}
We propose to use directly the decomposition on the reference basis regardless of the choice of hyperparameters; transferring the dependency on hyperparameters to the distribution of the coordinates. 
The inference problem has then the following general formulation,
\begin{equation}
 \pi_\mathrm{post}\left( \bs{\xi}, \bq \left| \bs{d}_\mathrm{obs}\right. \right) \propto \mc{L}(\bs{d}_\mathrm{obs} | \bs{\xi})\pi(\bs{\xi}|\bq)\pi_\mathbb{H}(\bq).
\end{equation}
This method therefore generalizes the handling of the autocovariance scaling in~\cite{marzouk2009} to the other hyperparameters. Note that, unlike the \gls{coc} method, the likelihood is independent of $\bq$. 
The advantage is to work with the reference coordinates $\bxi$ only, without using the $\bq$-dependent \gls{kl} decomposition,
\begin{equation}
 \ol{G}^r(\bx,\bxi) = \suml_{i=1}^{r} \ls{i}{r}\ol{u}_i(\bx)\xi_i, \text{ where } \bxi \sim \pi(\bxi|\bq).
\end{equation}
In this section, we derive the distribution for the $\bxi \coloneqq \left( \xi_i\right)_{1\leq i \leq r} \in \Xi \subset \R^r$ coordinates in order to best represent $G(\bx,\theta)$. The statistics of the new coordinates $\bxi$ are deduced by projecting the field on the reference basis. Assuming that $\ol{G}^r(\bx,\bxi) \simeq G(\bx,\theta)$, then
${\scal{ \ol{G}^r(\cdot,\bxi)}{\ol{u}_i}_\Omega} = {\scal{G(\cdot,\theta)}{\ol{u}_i}_\Omega}$
and using Eq.~\eqref{eq:KL-orthonormality} we obtain
\begin{equation}
 \ls{i}{r}\xi_i = \suml_{k=1}^{+\infty}\ls{k}{\bq}\scal{u_k(\cdot,\bq)}{\ol{u}_i}_\Omega \eta_k(\theta),
\end{equation}
where the coordinates $(\eta_k)_{k \in \N^\ast}$ are standard normal variables. By construction, $\bxi$ is a multivariate Gaussian random vector. It is clear that the $\xi_i$ have zero mean, 
\begin{equation}
 \E_\Theta( \xi_i) \simeq \ls[-]{i}{r}\suml_{k=1}^{+\infty}\ls{k}{\bq}\scal{u_k(\cdot,\bq)}{\ol{u}_i}_\Omega \E_\Theta(\eta_k) = 0.
\end{equation}
The covariance between $\xi_i$ and $\xi_j$ is 
\begin{align}
\E_\Theta(\xi_i \xi_j) &\simeq \E_\Theta\left( \left( \ls[-]{i}{r}\suml_{k=1}^{+\infty}\ls{k}{\bq}\scal{u_k(\cdot,\bq)}{\ol{u}_i}_\Omega \eta_k \right) \left( \ls[-]{j}{r}\suml_{k=1}^{+\infty}\ls{k}{\bq}\scal{u_k(\cdot,\bq)}{\ol{u}_j}_\Omega \eta_k\right)\right), \notag \\
&= \left(\ol{\lambda}_i \ol{\lambda}_j \right)^{-1/2}\suml_{k,k'=1}^{+\infty} \scal{u_k(\cdot,\bq)}{\ol{u}_i}_\Omega\scal{u_{k'}(\cdot,\bq)}{\ol{u}_j}_\Omega\ls{k}{\bq}\ls{k'}{\bq}\E_\Theta\left(\eta_k \eta_{k'}\right), \notag \\
&= \left(\ol{\lambda}_i \ol{\lambda}_j \right)^{-1/2} \scal{ \scal{ k(\cdot,\cdot,\bq)}{\ol{u}_i}_\Omega}{\ol{u}_j}_\Omega,
\end{align}
where the last line is obtained thanks to Mercer's theorem~\cite{mercer1909}. The reference coordinates $\bxi$ follow a prior normal distribution $\mc{N}(0,\Sigma(\bq))$, where the covariance matrix $\Sigma(\bq) \in \R^{r\times r}$ is defined by
\begin{equation}\label{eq:sigma-def}
\forall 1\leq i,j\leq r,\ \forall \bq \in \mathbb{H}, \qquad \Sigma(\bq)_{ij} = (\ol{\lambda}_i\ol{\lambda}_j)^{-1/2} \scal{\scal{k(\cdot, \cdot, \bq)}{\ol{u}_j}_\Omega}{\ol{u}_i}_\Omega.
\end{equation}
The covariance of the coordinates $\bxi$ for a given $\bq$ is obtained by projecting the $\bq$-dependent autocovariance on the span of the $r$-dimensional reference basis. This covariance matrix is invertible if $k(\bq)$ is positive definite, which is the case in the following applications. However, note that the conditionning of $\Sigma(\bq)$ can be quite bad if $r$ is very large. 

\begin{figure}[!ht]
  \centering
  \includegraphics[width=0.45\textwidth]{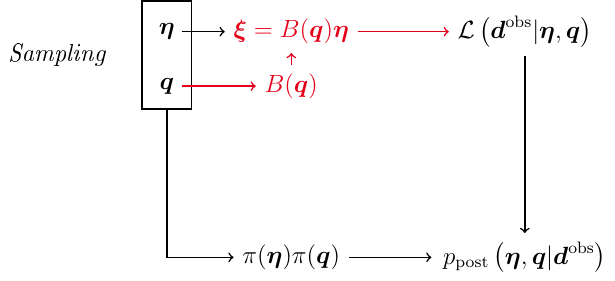}
  \includegraphics[width=0.45\textwidth]{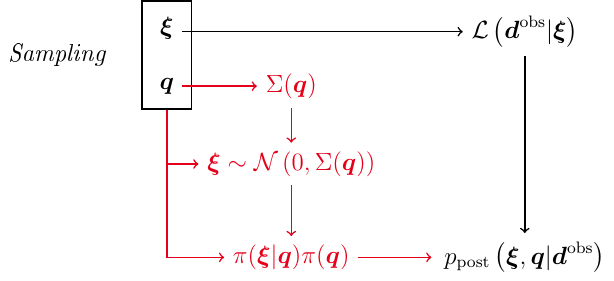}
  \caption{Workflow of the \gls{coc} (left) and the \gls{com} (right) methods. The quantities that are specific to the method are indicated in red.}\label{fig:workflow-methods}
\end{figure}

The principles of the \gls{coc} and \gls{com} methods are schematically represented on Fig.~\ref{fig:workflow-methods}.
For both methods, the model predictions only depends on the coordinates $\bxi$ in the reference basis.
The main advantage of the \gls{com} method is to split the posterior construction in two distincts parts: 
\begin{inparaenum}[i)]
\item the computation of the likelihood $\mc{L}(\bs{d}^\mathrm{obs}|\bxi)$ which is independent from the hyperparameters $\bq$ and 
\item the computation of the prior probability distribution $\pi(\bxi|\bq)\pi_\mathbb{H}(\bq)$ with the covariance matrix $\Sigma(\bq)$.
\end{inparaenum}
As the covariance matrix $\Sigma(\bq)$ only depends on the reference eigenelements which are computed offline, its computation is less expensive than the computation of $B(\bq)$ which requires to solve an eigenvalue problem at each step. In addition, the approximation of $\Sigma(\bq)$ by a surrogate model benefits from the regularity of the $\bq$-dependent autocovariance function. 
The \gls{com} representation also provides a better physical interpretation in comparison with the \gls{coc} one.
Indeed, the \gls{com} approach better distinguishes the coordinates from the hyperparameters, in the sense that changing the hyperparameters does not modify the field realization but its probability. 

\subsection{\gls{mcmc} for the change of measure}\label{subsec:com-implementation}
The logarithm of the posterior probability distribution is written as 
\begin{equation}
\log \pi_\mathrm{post}(\bxi,\bq | \bs{d}_\mathrm{obs}) \propto \log \mc{L}(\bs{d}_\mathrm{obs} | \bxi) + \log \pi(\bxi|\bq) + \log \pi_\mathbb{H}(\bq).
\end{equation}
A Gaussian likelihood is considered, assuming a centered \gls{iid} observation noise, with variance $\sigma_\varepsilon^2$: the noise level $\sigma_{\varepsilon}$ is the standard deviation of the observation noise $\varepsilon = \bs{d}_\mathrm{obs} - M(f_\mathrm{true})$ and is assumed to be \gls{iid} Gaussian,
\begin{equation}
\varepsilon \sim \mc{N}(0, \sigma_{\varepsilon}^2\mathrm{I}_{N}),
\end{equation}
where $\mathrm{I}_N$ is the identity matrix of size $N$. The likelihood writes
\begin{equation}\label{eq:likelihood}
\mc{L}(\bs{d}_\mathrm{obs}|\bxi) \coloneqq \inv{\sqrt{(2\pi\sigma_{\varepsilon}^2)^N}}\exp\left(-\frac{1}{2\sigma_{\varepsilon}^2} \norm{\bs{d}_\mathrm{obs} - M(\bxi)}_{l_2}^2\right).
\end{equation}
Other likelihood functions and error models can be used~\cite{zhang2014,leoni2022}. 
Since $\bxi \sim \mc{N}(0,\Sigma(\bq))$, the prior log-density of the random vector $\bxi$ is
\begin{equation}
\log \pi(\bxi|\bq) \propto -\inv{2} \left( \mathrm{logdet}_\Sigma(\bq) + \bxi(\bq)^\top\Sigma(\bq)^{-1}\bxi(\bq) \right),
\end{equation}
where $\mathrm{logdet}_\Sigma(\bq) \coloneqq \mathrm{log}\left(\mathrm{det}\left[ \Sigma(\bq) \right] \right)$.
The noise level $\sigma_{\varepsilon}$ is estimated during the inference, assuming a Jeffreys prior~\cite{jeffreys1946,kass1996a,chopin2009}, $\pi_{\sigma_{\varepsilon}}(\sigma_{\varepsilon}) \propto 1/\sigma_{\varepsilon}$.

The \gls{com} posterior has a hierarchical structure~\cite{betancourt2013}. In particular, since $\pi(\bxi|\bq)$ can be highly sensitive to $\bq$, it can be difficult to define a proposal structure suitable over the whole sampling space. The \gls{mh} algorithm, even in its adaptive version~\cite{haario2001, roberts2009}, can be challenged by the \gls{com} hierarchical structure without a proposal distribution $\pi_\mathrm{tr}$ depending on the current state of the chain $(\bxi^{(n)},\bq^{(n)})$. Our strategy to overcome this difficulty consists in sampling an auxiliary variable $\ol{\bxi}$ independent from $\bq$, using a fixed proposal structure and a $\bq$-dependent linear transformation to obtain $\bxi$. Specifically, we consider $\ol{\bxi} \sim \mc{N}(0,\Sigma_{\ol{\bxi}})$ and set $\bxi$ through 
\begin{equation}
  \bxi = \Sigma(\bq)^{1/2}\Sigma_{\ol{\bxi}}^{-1/2}\ol{\bxi}.
\end{equation}

At each \gls{mcmc} step, the \gls{mh} acceptance criterion is used. Denoting $\bs{Y}^{(n)} = (\bxi^{(n)},\bq^{(n)}, \sigma_\varepsilon^{(n)})$ the current state and $\bs{Y}^\ast = (\bxi^\ast,\bq^\ast,\sigma_\varepsilon^\ast)$ the proposed one, the proposition is accepted with probability 
\begin{equation}
p_\mathrm{MH} = \mathrm{min}\left(\frac{ \pi_\mathrm{post}(\bs{Y}^\ast)\pi_\mathrm{tr}(\bs{Y}^{(n)}|\bs{Y}^\ast)}{\pi_\mathrm{post}(\bs{Y}^{(n)})\pi_\mathrm{tr}(\bs{Y}^\ast|\bs{Y}^{(n)}) }, 1\right),
\end{equation}
where $\pi_\mathrm{tr}(\bs{Y}|\bs{X})$ is the transition probability from $\bs{X}$ to $\bs{Y}$. In our case, this transition probability is not symmetric,
\begin{equation}
  \frac{\pi_\mathrm{tr}(\bs{Y}^{(n)}|\bs{Y}^\ast)}{\pi_\mathrm{tr}(\bs{Y}^\ast|\bs{Y}^{(n)})} = \left( \frac{\det{\Sigma(\bq^\ast)}}{\det{\Sigma(\bq^{(n)})}}\right)^{1/2}.
\end{equation}
Note that $\det{\Sigma(\bq)}$ is computed for the prior, such that the evaluation of the transition probability ratio add no computational cost. Further, as $\Sigma(\bq)$, $\Sigma(\bq)^{1/2}$ is independent of the indexation and orientation of the eigenelements.

\section{Accelerating sampling in the \gls{com} framework}\label{sec:acceleration}
The Bayesian inference of $\ol{G}^r(\bx,\bxi)$ with $\bs{\xi}\sim \mc{N}(0,\Sigma(\bq))$ requires the sampling of the $r+|q|$-dimensional space of $(\bxi,\bq) \in \Xi\times\mathbb{H}$ with a \gls{mcmc} algorithm.
As previously mentioned, one advantage of the \gls{com} is that the predictions, and therefore the likelihood, only depends on the coordinates $\bxi$. In this section, we briefly present \gls{pc} expansions used to replace the set of forward model predictions. The coefficients and determinant of the \gls{com} matrix used to define the likelihood function are also approximated with \gls{pc} expansions.

\subsection{Polynomial chaos surrogates}
Surrogate models are widely used in uncertainty propagation and Bayesian inference~\cite{peherstorfer2018,al-ghosoun2021}. We rely on \gls{pc} expansions~\cite{wiener1938a,ghanem1991a,xiu2002} which have been used in various fields~\cite{marzouk2009,laloy2013,sochala2021a,han2022,robbe2023,khatoon2023,meles2024}. Let $f$ be a second order functional of a $r$-dimensional vector $\bz$ with independent components such that $p_{\bz}(\bs{y}) = \suml_{i=1}^r p_{i}(y_i)$. 
The \gls{pc} expansion $\widetilde{f}$ of $f$ is a linear combination of polynomials in $\bz$,
\begin{equation}\label{eq:pc-expansion}
f(\bz) \simeq \widetilde{f}(\bz) = \suml_{a \in \mc{A}}f_a \Psi_a(\bz),
\end{equation}
where the \gls{pc} basis functions $\{\Psi_a(\bz)\}_{a\in\mc{A}}$ are multivariate orthogonal polynomials with respect to the density of $\bz$, $\{f_a\}_{a\in \mc{A}}$ are deterministic \gls{pc} coefficients and $\mc{A}$ is a set of multi-indexes. The \gls{pc} basis functions are the product of orthonormal univariate polynomials, $\Psi_{a = (a_1, \dots a_d)}(\bz) = \Pi_{i}\Psi_{a_i}^i(\bz_i)$, where $\Psi_{a_i}^i$ is a polynomial of degree $a_i$ in $\bz_i$.
Several methods can be implemented to estimate the \gls{pc} coefficients. In this work, we use the \gls{psp} method to approximate the projection coefficients $\{f_a\}_{a\in \mc{A}}$ of $f(\bz)$~\cite{orszag1972,constantine2012}. The \gls{psp} method relies on sequences of nested $1$-dimensional discrete projection operators into spaces of increasing polynomial degrees. The $1$-d projections are sparsely tensorized, using the Smolyak formula. The \gls{psp} yields sparse projection operators, free of internal aliasing for $f \in \mathrm{Span}\{\Psi_a, \ a \in |\mc{A}|\}$, where $\mc{A}$ depends on the composing sequences and their retained tensorizations. In practice, the \gls{psp} method uses evaluations of $f$ at the nodes $\bz^{(l)}$ of a sparse grid $S = \{\bz^{(k)}\}_{1\leq k \leq N_\mathrm{PSP}}$ and computes the \gls{pc} coefficients $f_a$ through 
\begin{equation}
   \forall a \in \mc{A}, \quad f_a = \suml_{k=1}^{N_\mathrm{PSP}}\Pi_{ak}f(\bz^{(k)}).
\end{equation}

The surrogate model accuracy is assessed by estimating the \gls{rrmse} with a validation set of $N_v$ samples $\mc{Z} = \left\{\bz^{(i)}\right\}_{1\leq i \leq N_v}$ drawn from $p_{\bz}$, 
\begin{equation}\label{eq:def-rrmse}
\mathrm{RRMSE}(f,\widetilde{f}) = \sqrt{ \suml_{\bz \in \mc{Z}} \norm{f(\bz) -\widetilde{f}(\bz)}^2 / \suml_{\bz \in \mc{Z}} \norm{f(\bz)}^2}.
\end{equation}

\subsection{Surrogate models for likelihood computation}
As seen in Sections~\ref{subsec:com} and~\ref{subsec:com-implementation}, the model predictions in the likelihood only depend on the coordinates of the field $\bxi$. Therefore, we construct a surrogate model of the prediction vector $\bs{d}(\bxi)$ with the \gls{psp} method: 
\begin{equation}
   \bs{d}(\bxi) \simeq \suml_{a \in \mc{A}}d_a \Psi_a(\bxi).
\end{equation}
Recall that the coordinates $\bxi$ follow $\mc{N}(0, \Sigma(\bq))$ in the \gls{com} method. To avoid the $\bq$-dependency of $\bxi$, we consider that $\bxi\sim\mc{N}(0,\mathrm{I}_r)$ for the \gls{pc} surrogate construction. This choice corresponds to a simple approximation of the marginal coordinates distribution over the hyperparametric domain since $\E_{\mathbb{H}}\left(\E_\Theta(\bxi\left|\bq\right.)\right) = 0$ and $\E_{\mathbb{H}}(\Sigma(\cdot)) = \mathrm{I}_r$ from equations~\eqref{eq:Bref-def},\eqref{eq:sigma-def}. Note that this is not the true marginal coordinates distribution, since the sum of Gaussian distributions is not Gaussian ( see~\ref{appendix:xi-law} for more details). The model predictions at the \gls{psp} points are computed by
\begin{inparaenum}[i)]
 \item building the field associated with the \gls{psp} point $g^{(l)}(\bx) = \suml_{i=1}^r \ol{\lambda}_i^{1/2}u_i(\bx)\xi_i^{(l)}$,
 \item solving the forward model with $g^{(l)}$ and evaluate $\bs{d}^{(l)} = d(\bxi^{(l)})$.
\end{inparaenum}

\subsection{Surrogate models for prior computation}\label{subsec:acceleration:surrogate-prior}
The \gls{com} formulation and its sampling require the computation of $\Sigma(\bq)^{-1}$, $\mathrm{logdet}_\Sigma(\bq)$, and $\Sigma(\bq)^{1/2}$ at each step. We construct surrogate models of these three quantities to accelerate the \gls{mcmc} sampling.
To ensure that the approximated $\Sigma(\bq)^{-1}$ is non negative, we approximate its square root. To ensure the positivity of the $\Sigma(\bq)^{-1}$ surrogate, we adopt the approach of~\cite{reis2021} for the \gls{pc} approximations of semi positive definite operators in the context of domain decomposition methods for stochastic partial differential equations. The approach is based on the spectral decomposition of $\Sigma(\bq)$ in 
 \begin{equation}
   \Sigma(\bq) = U(\bq)\Lambda(\bq)U(\bq)^\top, \text{ with } U(\bq)U(\bq)^\top = U(\bq)^\top U(\bq)=\mathrm{I}_r,
 \end{equation}
 and defines
 \begin{equation}
   \Sigma(\bq)^{\pm 1/2} \coloneqq U(\bq)\Lambda(\bq)^{\pm 1/2}U(\bq)^\top.
 \end{equation}
 The advantage of this approach is that the \gls{psp} approximation 
 \begin{equation}
  \widetilde{\Sigma}^{\pm 1/2} = \suml_{a \in \mc{A}} \Sigma_a^{\pm 1/2}\Psi_a(\bq), \text{ with } \Sigma_a^{\pm 1/2} = \suml_{k=1}^{N_\mathrm{PSP}} \Pi_{ak}\Sigma(\bq^{(k)})^{\pm 1/2}
 \end{equation}
 is not affected by the indexation and orientation of the $\Sigma(\bq^{(k)})$ eigenelements at the different \gls{sg} points. Finally, we use 
 \begin{equation}
   \widetilde{\Sigma}^{-1}(\bq) \coloneqq \widetilde{\Sigma}^{-1/2}\widetilde{\Sigma}^{-1/2}.
 \end{equation}
 In addition, the \gls{psp} approximation of $\mathrm{logdet}_\Sigma(\bq)$ writes as 
 \begin{equation}
   \logdetsurr_\Sigma(\bq) = \suml_{a \in \mc{A}}l_a \Psi_a(\bq), \text{ with } l_a = \suml_{k=1}^{N_\mathrm{PSP}}\Pi_{ak}\mathrm{logdet}_\Sigma(\bq^{(k)}).
 \end{equation}

When the prior range is large and the number of modes $r$ is high, the prior $\pi(\bxi|\bq)$ can become highly stretched with respect to $\bq$. In that case, the \gls{pc} approximation of the conditional law can be difficult and requires a high \gls{psp} level. In practice, we adapt the \gls{psp} method level to ensure a \gls{rrmse} of less than $0.1\%$ on the surrogates of $\Sigma(\bq)^{-1}$ and $\mathrm{logdet}_\Sigma(\bq)$. As for the surrogate predictions, the \gls{pc} error should be dominated by the observation error. This is verified~\textit{a posteriori}, and the analysis is led with a higher \gls{psp} level if needed.

\section{Application to a transient diffusion problem}\label{sec:transdiff}
In this section, the \gls{com} method is applied to infer a diffusivity field in the one-dimensional \gls{td} problem presented in~\cite{sraj2016}. This problem is summarized in Section~\ref{sec:td:setup-pb}, the behavior and convergence of the \gls{pc} surrogates of $\Sigma^{\pm 1/2}$ and $\mathrm{logdet}_\Sigma$ are analysed in Section~\ref{sec:td:surr-an} and the inference results are presented in Section~\ref{sec:td:inference-res}, where a comparison with the \gls{coc} method is proposed.

\subsection{Case presentation}\label{sec:td:setup-pb} 
The 1D \gls{td} equation writes 
\begin{equation}\label{eq:diffusion}
\forall t\in (0,T=0.05),\ \bx \in \Omega = (0,1),\quad \frac{\partial U(\bx,t)}{\partial t} = \frac{\partial}{\partial \bx}\left(\nu(\bx) \frac{\partial U(\bx,t)}{\partial \bx}\right),
\end{equation}
with the boundary and initial conditions,
\begin{equation}
U(\bx=0,t)=-1,\quad U(\bx=1,t)=1,\quad\text{and}\quad U(\bx,t=0)=0.
\end{equation}
The field $\nu$ is such that $0<\nu(\bx)<+\infty$ almost everywhere in $\Omega$, is unknown and to be inferred.
For a given $\nu$, the transient diffusion equation is solved with a $\mathbb{P}_1$-finite element method in space and a second order implicit time-integration scheme.
The objective of the inference is to learn the diffusitivity field $\nu$ from observations of $U$. The observations consist in $N=N_x\times N_t=18\times13=234$ noisy evaluations of $U$, uniformly distributed in space and time. They are obtained by solving Pb.~\eqref{eq:diffusion} and then synthetically corrupted by random independent Gaussian centered measurement noises $\varepsilon \sim \mc{N}(0, \sigma_{\varepsilon}^2\coloneqq0.01)$. For the inference, we consider \textit{a priori} a log-normal stationnary field $\nu$, and then infer $g\coloneqq \log \nu$ using a Gaussian prior with zero mean and with a squared exponential autocovariance function $k(\bq)$,
\begin{equation}
\forall (\bx,\by)\in \Omega^2,\quad\bq=\{A,l\}\in\mathbb{R}_+^2,\quad k\left(\bx,\by,\bq\right)\coloneqq A\exp\left(-\norm{\bx - \by}^2/(2l^2) \right),
\end{equation}
where $A$ denotes the amplitude and $l$ the correlation length. The prior distributions of these two hyperparameters are 
\begin{equation}
A \sim \mathrm{InvGamma}(3,1)\quad\text{and}\quad l\sim \mathrm{log}\text{-}\mc{U}(0.1,0.7).
\end{equation}
The decomposition in the reference basis shown in Eq.~\eqref{eq:ref-decomposition} leads to the following approximation, 
\begin{equation}\label{eq:td:proj-kl}
 g_\mathrm{true} \simeq g_\mathrm{kl}(\bxi) = \ol{G}^r(\cdot,\bxi), \text{ with } \xi_i = \scal{\ol{\lambda}_i^{-1/2}\ol{u}_i}{g_\mathrm{true}}_\Omega.
\end{equation}
The number of coordinates $r$ is set equal to $8$ and captures $99.8\%$ of the field prior variance.

\subsection{Behavior of the change of measure}\label{sec:td:surr-an}
The aim of this section is to analyze the behavior of the \gls{com}'s surrogates for the \gls{td} case. It focuses on the analysis of $\Sigma(\bq)$ and its derived quantities $\Sigma^{\pm 1/2}(\bq)$ and $\mathrm{logdet}_\Sigma(\bq)$ when $l$ varies, since the $A$-dependency is reduced to a multiplicative factor. 

\begin{figure}
  \centering
  \begin{minipage}[t]{0.3\textwidth}
      \centering
      \includegraphics[height=0.2\textheight]{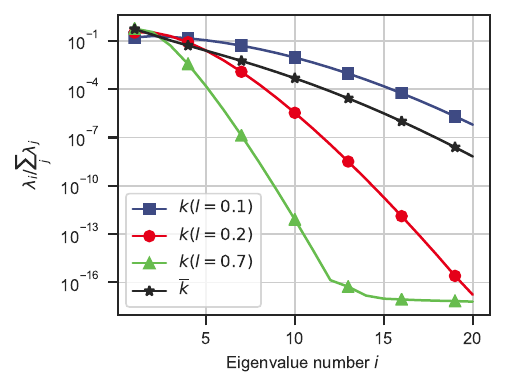}
      \caption{\acrshort{td} case - Eigenvalues decay for different correlation lengths $l$ as indicated.}\label{fig:eigenvalues}
  \end{minipage}\hfill
  \begin{minipage}[t]{0.3\textwidth}
      \centering
      \includegraphics[height=0.2\textheight]{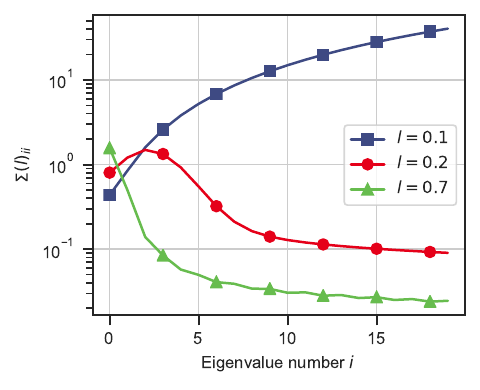}
      \caption{\acrshort{td} case - Diagonal terms of $\Sigma(l)$ corresponding to the prior variance of $\bxi$.}\label{fig:prior-xi}
  \end{minipage}\hfill 
  \begin{minipage}[t]{0.3\textwidth}
    \centering
    \includegraphics[height=0.2\textheight]{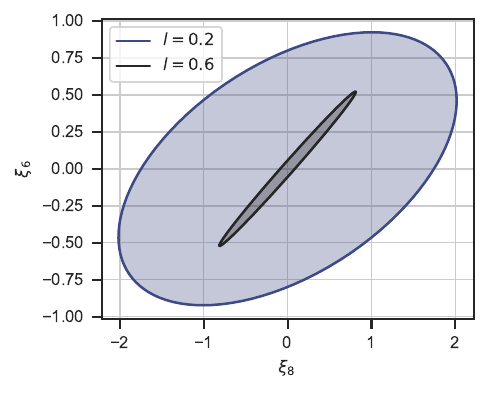}
    \caption{\acrshort{td} case - Prior Covariance $\Sigma(l)$ for two correlation lengths. Projection on the $6th$ and $8th$ coordinates. }\label{fig:sampling-hatcov-cov}
\end{minipage}
\end{figure}

\begin{figure}
  \centering
  \subfloat{
    \includegraphics[width=0.35\textwidth]{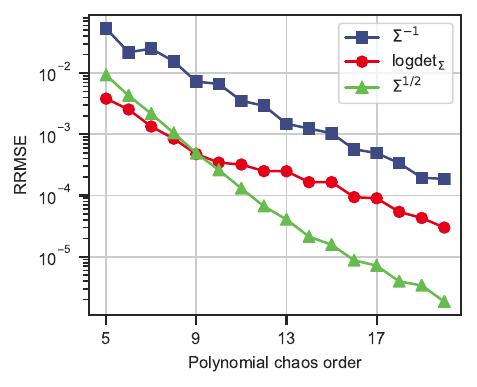}
  } \hspace{0.2cm}
  \subfloat{
    \includegraphics[width=0.35\textwidth]{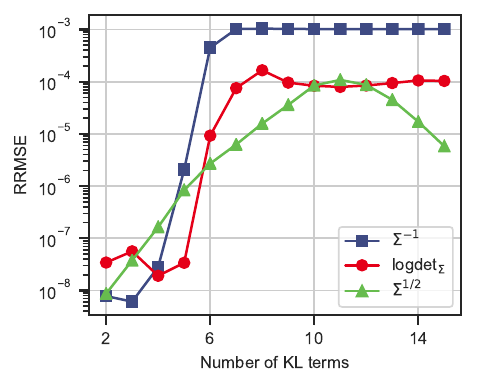}
  }
  \caption{\acrshort{td} case - \gls{rrmse} values (left) according \gls{pc} order ($r=8$), (right) according number of terms (\gls{pc} order $=15$). $1,000$ random validation points.}\label{fig:sigma-cv-rrmse} 
\end{figure}

The eigenvalues of the covariance for different correlation lengths are plotted on Fig.~\ref{fig:eigenvalues}: the higher the correlation length, the faster the decay. In particular, for large correlation lengths, the information of the $l$-dependent basis is concentrated in the few first modes. The spectrum of the averaged covariance (\textit{i.e.} leading to the reference basis) is intermediate. Figures~\ref{fig:prior-xi} and~\ref{fig:sampling-hatcov-cov} show how changes in $l$ affect the prior coordinates $\bxi$. Figure~\ref{fig:prior-xi} reports the diagonal coefficients $\Sigma_{ii}(\bq)$ for three values of $l$:
the higher the correlation length, the smaller the variance of the last coordinates. A coordinate variance is close to zero when the associated mode is irrelevant for the approximation of the prior field. 
Figure~\ref{fig:sampling-hatcov-cov} shows the resulting $90\%$ confidence interval for $\xi_6$ and $\xi_8$ and two values of $l$, highlighting the strong dependency of $\Sigma(\bq)$ on $l$. The \gls{com} induces a high sensitivity of the $\bxi$ prior distribution to $l$.
Tab.~\ref{tab:sigmaanalysis} highlights the influence of the prior range and the number of coordinates on the eigenvalues and covariance matrix of the \gls{com} method. Two conclusions can be drawn: 
\begin{inparaenum}[i)]
\item diminishing $r$ increases the error on the field prior, and diminishing $l_-$ builds a reference basis that is closer to small correlation length bases, each effect increasing the truncation error and 
\item increasing the distance between $l_-$ and $l_+$ or the number of terms $r$ yields a higher $l^2$-norm.
\end{inparaenum}
The numerical invertibility of the covariance matrix is compromised when the hyperparameter space or the number of coordinates is too large.
\begin{table}
  \centering
\begin{tabular}{|cccccc|}\hline 
  $r$ & $l_-$ & $l_+$ & $\sum_{i=1}^r\ol{\lambda}_i/\sum_{i=1}^{+\infty}\ol{\lambda}_i$ & $\norm{\Sigma(\bq)^{-1}}_{l_+^2}$ & \gls{rrmse}($\Sigma(l)^{-1}$) \\ \hline
  $8$ & $0.1$ & $0.7$ & $99.8\%$ & $6.7{\times}10^3$ & $1.0{\times}10^{-3}$ \\ 
  $8$ & $0.05$ & $0.7$ & $97.1\%$ & $3.3{\times}10^4$ & $6.4{\times}10^{-3}$ \\ 
  $8$ & $0.1$ & $0.9$ & $99.8\%$ & $1.1{\times}10^4$ & $1.5{\times}10^{-3}$ \\ 
  $6$ & $0.1$ & $0.7$ & $99.0\%$ & $6.4{\times}10^3$ & $4.5{\times}10^{-4}$ \\ \hline
\end{tabular}
\caption{\acrshort{td} case - Quantities of interest for different numbers of coordinates and prior ranges. \gls{rrmse} is computed for \gls{pc} order equal to $15$.}\label{tab:sigmaanalysis}
\end{table}
Finally, the \gls{rrmse} of $\Sigma^{1/2},\Sigma^{-1}$ and $\mathrm{logdet_{\Sigma}}$ are plotted on Fig.~\ref{fig:sigma-cv-rrmse} for different \gls{pc} order and number of terms. For $r=8$, a \gls{pc} order of $15$ ensures that the error is below $0.1\%$ for the three quantities. A Tikhonov regularization is implemented in order to prevent numerical instability when increasing the number of terms.

\subsection{Inference results}\label{sec:td:inference-res}
\begin{figure}
  \centering 
  \subfloat{
    \includegraphics[width=0.45\textwidth]{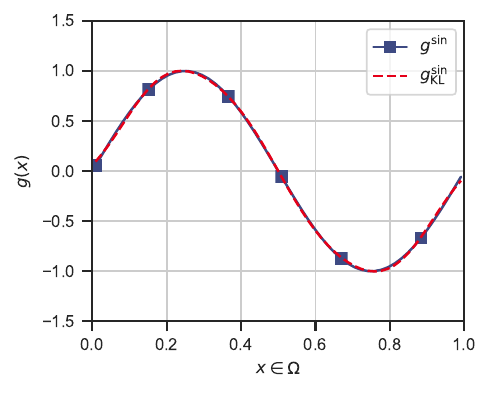}
  }\hfill\subfloat{
    \includegraphics[width=0.45\textwidth]{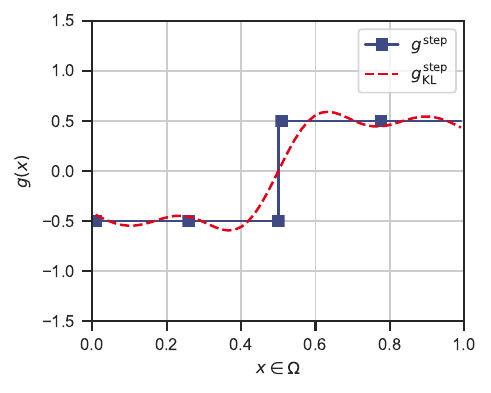}
  }
  \caption{\acrshort{td} case - True diffusivity fields $g_\mathrm{true}$ and their projections $g_\mathrm{kl}$. (left) $g^\mathrm{sin}$, (right) $g^\mathrm{step}$ }\label{fig:transdiff:bestproj}
\end{figure}

We test the \gls{com} for two different true log-diffusivity fields, namely
\begin{inparaenum}[i)]
\item a sinusoidal profile: $g^\mathrm{sin}(\bx) = \mathrm{sin}(2\pi\bx)$ and
\item a step function: $g^\mathrm{step}(\bx) = \left\{\begin{array}{l} -1/2, \ \text{if } \bx <0.5, \\ 1/2 \ \text{else}\end{array}\right.$
\end{inparaenum}.
The projections of the true fields in the reference basis (Eq.~\eqref{eq:td:proj-kl}) are shown on Fig.~\ref{fig:transdiff:bestproj}. As expected, the field $g^\mathrm{sin}$ is well approximated in the reference basis whereas the field $g^\mathrm{step}$ is considerably smoothed and presents oscillations, because this non-smooth field is not a likely realization from the prior.

\begin{figure}
  \centering 
  \includegraphics[width=0.25\textwidth]{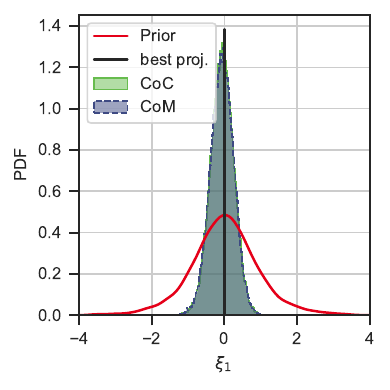}
  \includegraphics[width=0.25\textwidth]{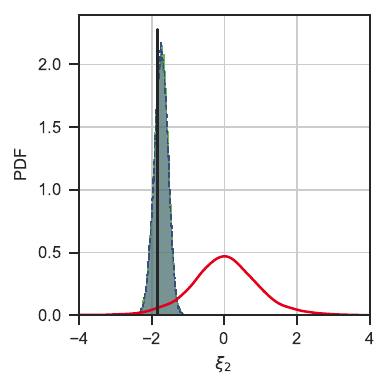}
  \includegraphics[width=0.25\textwidth]{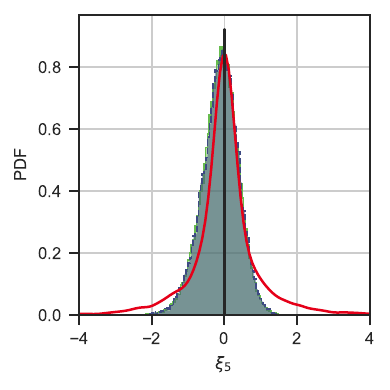}\\
  \includegraphics[width=0.25\textwidth]{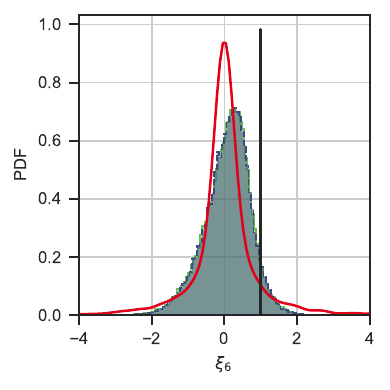}
  \includegraphics[width=0.25\textwidth]{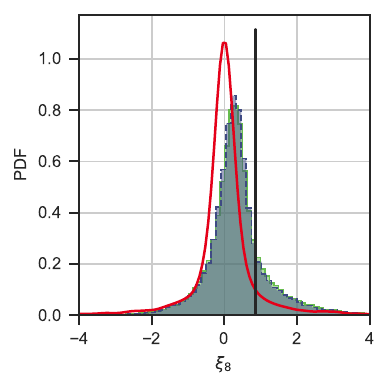}
  \includegraphics[width=0.25\textwidth]{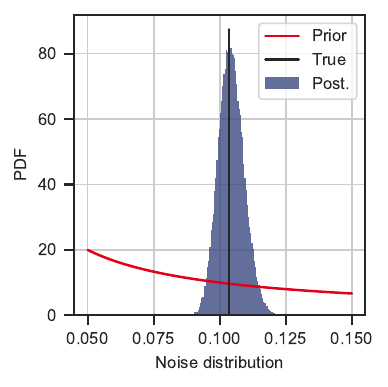}
  \caption{\acrshort{td} case for $g^\mathrm{sin}$ - Marginal posterior distribution of coordinates $1$, $2$, $5$, $6$ and $8$ and noise level  (coordinates $3$, $4$ and $7$ are not shown for brevity). Best projection corresponds to the coordinates used to plot Fig.~\ref{fig:transdiff:bestproj} and true noise is the standard deviation of the synthetical noises.}\label{fig:transdiff:postkl}
\end{figure}

We use the \gls{mcmc} presented in Section~\ref{subsec:com-implementation} to sample the posterior distribution. The \gls{coc} method introduced in Section~\ref{subsec:com-coc} is also implemented for validation purposes. The proposal matrix is adapted each $2.5{\times}10^4$ steps during a burning phase of $2.5{\times}10^5$ steps, then $1{\times}10^6$ steps are used to sample the posterior distribution. This choice leads to a multi effective sample size~\cite{vats2019} around $10,000$ for both \gls{com} and \gls{coc} methods. We start the analysis by examining the field coordinates. The marginal posterior distributions of five coordinates are plotted on Fig.~\ref{fig:transdiff:postkl} for the $g^\mathrm{sin}$ field. 
For the comparison, the \gls{coc} samples $(\bs{\eta},\bq)$ are translated into samples of $\bxi$ using the \gls{coc} $\bxi = B(\bq)\bs{\eta}$. The $\bxi$ marginal prior distribution $\pi(\bxi) = \Int_{\mathbb{H}}\pi(\bxi|\bq)\pi_\mathbb{H}(\bq)d\bq$ is also plotted as reference. 
We observe that the two methods yield similar results, with the posterior distributions for the first $5$ coordinates centered on the best projections of the true field. The other coordinates are less informed by the observations. The noise level posterior distribution is peaked around the true value.

\begin{figure}
  \centering 
  \subfloat{\includegraphics[width=0.45\textwidth]{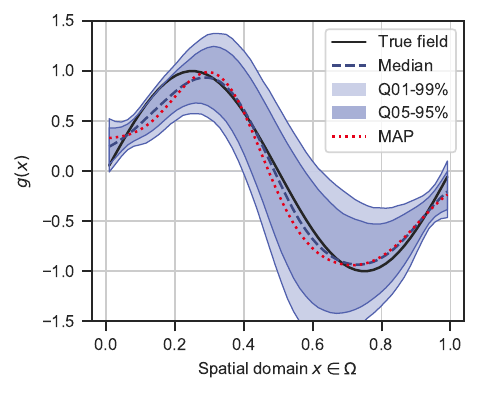}}
  \hfill
  \subfloat{\includegraphics[width=0.45\textwidth]{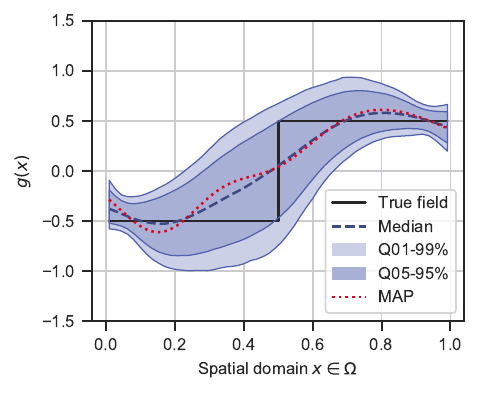}}
  \\\vspace{-0.5cm}
  \subfloat{\includegraphics[width=0.45\textwidth]{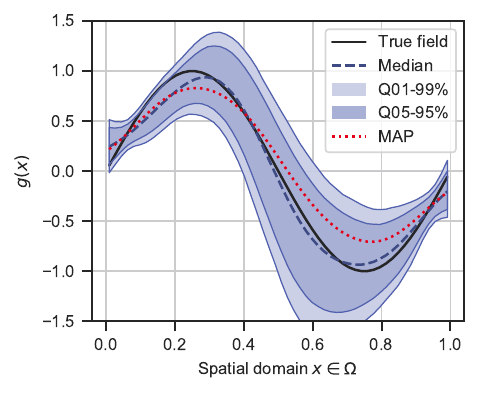}}\hfill
  \subfloat{\includegraphics[width=0.45\textwidth]{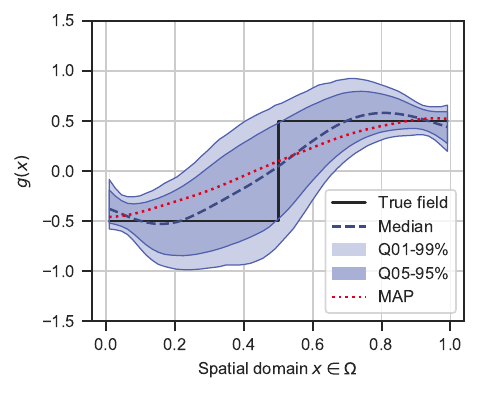}}
  \caption{\acrshort{td} case - Posterior distribution of the inferred fields. Black line: $g_\mathrm{true}$; star blue line: median field; dashed red line: \gls{map}; blue intervals: $1\%-99\%$ and $5\%-95\%$ quantiles. (left) $g^\mathrm{sin}$ case / (right) $g^\mathrm{step}$ case ; (top) Change of coordinates / (bottom) Change of measure}\label{fig:transdiff:fielddistrib}
\end{figure}

\begin{figure}
  \centering
  \subfloat{
  \includegraphics[width=0.3\textwidth]{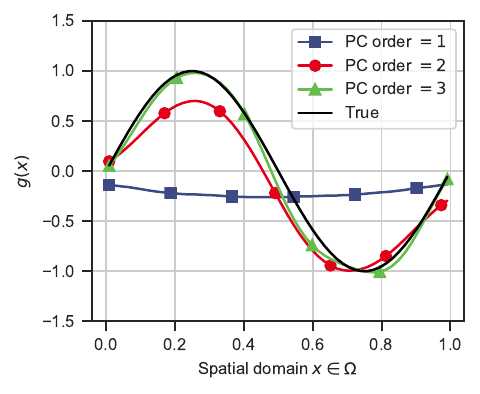}}\hfill
  \subfloat{\includegraphics[width=0.3\textwidth]{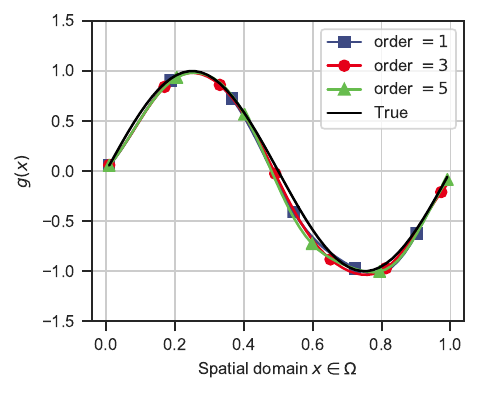}}\hfill
  \subfloat{ \includegraphics[width=0.3\textwidth]{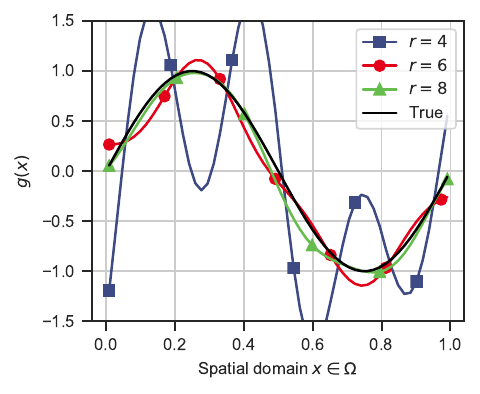}}  
  \caption{\acrshort{td} case for $g^\mathrm{sin}$ - Convergence of \gls{mcmc} posterior median according to different parameters (\gls{com} method): (left) \gls{pc} order for the forward model, (middle) \gls{pc} order for the \gls{com} prior quantities, (right) number of coordinates. For fixed parameters, inference realized with $r=8$, $\sigma_{\varepsilon}=0.02$, \gls{pc} order $\widetilde{\mc{M}}=5$, \gls{pc} order prior $=15$.}\label{fig:transdiff:cplt-analysis}
\end{figure}

Figure~\ref{fig:transdiff:fielddistrib} presents the posterior distributions confidence intervals of the inferred fields. 
The mean, median, \gls{map} as well as the $1\%$-$99\%$ and $5\%$-$95\%$ quantiles are plotted. The \gls{map} corresponds here to the chain sample maximizing the posterior distribution.
As for the coordinates estimations, the results of the two methods are very close for $g^\mathrm{sin}$ and $g^\mathrm{step}$. 
Note that the $g^\mathrm{sin}$ medians are shifted relatively to the true field, this gap is reduced when decreasing the noise level (see Fig~\ref{fig:transdiff:cplt-analysis} where $\sigma_\varepsilon$ is set to $0.02$).
Figure~\ref{fig:transdiff:cplt-analysis} assesses the convergence of the \gls{com} inference result with respect to the \gls{pc} order for the forward model, the number of coordinates and the \gls{pc} order for the \gls{com} prior quantities. 
Capturing the stochastic nonlinearities of the forward model is crucial (left figure) while the median field appears to be less sensitive to the approximation of the \gls{com} prior quantities (middle figure). 
It is clear that a too low number of coordinates is detrimental to the field inference (right figure).

\section{Application to seismic tomography}\label{sec:tomo}
The \gls{com} method is applied to a \acrfull{st} problem based on the propagation of the P-waves in the first kilometer of the Earth's crust. 
After introducing the test case in Section~\ref{sec:tomo-presentation}, the results are presented in Section~\ref{sec:tomo-online}. Results for different field shapes are presented in Section~\ref{sec:tomo:comp-field}. Section~\ref{sec:comp-param} provides a comparison with inference results with fixed hyperparameters.

\subsection{Case presentation}\label{sec:tomo-presentation}
We consider the inference of a $2$D continuous seismic velocity field $v$ from observations of traveltimes between known sources and receivers. The observation vector $\bs{d}_\mathrm{obs}$ corresponds to the $5$ sources and $23$ receivers depicted in Fig.~\ref{fig:mcmc-field}, such that the dimension of $\bs{d}_\mathrm{obs}$ is $N=115$.
The forward model is the eikonal equation~\cite{noble2014} that relates the traveltimes map $t_s(\bx)$ from a source at $\bx = s$ to the velocity field $v$,
\begin{equation}\label{eq:eikonal}
\lvert \nabla t_s(\bx)\rvert^2 = \inv{v^2(\bx)}, \text{ with } t_s(s) = 0
\end{equation}
where $\bx\in\Omega $ is the spatial position. The quantities are expressed with SI units: $t\sim [\mathrm{s}]$, $x\sim[\mathrm{m}]$ and $v\sim[\mathrm{m}\cdot \mathrm{s}^{-1}]$.
The observations $\bs{d}_\mathrm{obs}$ are synthetically generated by solving Eq.~\eqref{eq:eikonal} with $v=v_{\mathrm{true}}$ and corrupting the traveltimes with an \gls{iid} additive noise having a centered Gaussian distribution with standard deviation equal to $\sigma_{\varepsilon} = 0.002\sim[\mathrm{s}]$,
\begin{equation}
  \bs{d}_\mathrm{obs} = \bs{d}_\mathrm{true} + \varepsilon, \quad \varepsilon = \mc{N}(0, \sigma_\varepsilon^2\mathrm{I}_N).
\end{equation}
For the sources-receivers geometry and the true field we consider, the noise is around $1\%$ of the true times. 
In this work, $v_\mathrm{true}:\Omega \rightarrow \R_+$ is adapted from~\cite{obrien1994} and depends only on the vertical coordinate. We again consider a log-normal prior distribution:
\begin{equation}
  v(\bx) = \exp(g(\bx)), \quad g \sim \mc{N}(c,k), \quad k(\bx,\bs{y}) = A\exp\left(\frac{-\norm{\bx - \bs{y}}^2}{2l^2} \right),
\end{equation}
where $c \in \R$ is a constant trend to infer. The field $g(\bx)$ writes 
\begin{equation}
  g(\bx) = c+\suml_{i=1}^r \ol{\lambda}_i^{1/2} \ol{u}_i(\bx)\xi_i,
\end{equation}
with $(\ol{u}_i,\ol{\lambda}_i)_{1\leq i\leq r}$ the eigenelements of $\E_\mathbb{H}(k)$. In that case, the forward model surrogate depends on the coordinates $\bxi$ and on the constant trend $c$. \textit{A priori}, the constant trend is independent from the coordinates and from the other hyperparameters, leading to the following expression for the posterior distribution 
\begin{align*}
  \pi_\mathrm{post}(\bxi, \bq, c, \sigma_\varepsilon) &\propto \mc{L}(\bs{d}^\mathrm{obs}|\bxi,\bq,c, \sigma_\varepsilon)\pi(\bxi|\bq,c,\sigma_\varepsilon)\pi(\bq,c,\sigma_\varepsilon) \\
  &\propto \inv{\sqrt{(2\pi\sigma_\varepsilon^2)^N}}\exp\left( \frac{-\norm{\bs{d}^\mathrm{obs} - \bs{d}(\bxi,c)}_{l_2}^2}{2\sigma_\varepsilon^2}\right) \times \inv{\sqrt{(2\pi)^r\left|\Sigma(\bq)\right|}}\exp\left(\frac{-\bxi^\top \Sigma(\bq)^{-1}\bxi}{2}  \right) \\
  &\quad \times \pi(c)\pi_\mathbb{H}(\bq)\pi_{\sigma_{\varepsilon}}(\sigma_\varepsilon)
\end{align*} 
The prior distributions of the constant trend and the two hyperparameters are 
\begin{equation}
c\sim\mc{U}(6.9,8.1),\quad A \sim \mathrm{InvGamma}(21,1)\quad\text{and}\quad l\sim\mc{U}(10,100).
\end{equation}
These priors have been chosen to generate a large set of plausible velocity fields, see Fig.~\ref{fig:mcmc-field} (left). The number of reduced coordinates is $r=20$; it accounts for $96.8\%$ of the prior variance and allows building accurate and tractable surrogate models with sparse grids. 
For $r=20$, the best approximation of $v_\mathrm{true}$, 
\begin{align*}
  v_\mathrm{best}^r(\bx) &= \exp \left( c_b + \suml_{i=1}^r \ol{\lambda}_i^{1/2}\ol{u}_i(\bx)v_i \right), 
\end{align*}
where $c_b$ is the spatial mean of $\log v_\mathrm{true}$ and the coordinates $v_i$ result from the projection of $\log v_\mathrm{true} - c_b$ on the reference basis as done in Eq.~\eqref{eq:ref-decomposition}. This approximation
yields to low relative mean squared errors
\begin{equation}
  \norm{v_\mathrm{true}-v_\mathrm{best}^r}/\norm{v_\mathrm{true}} \simeq 2\% \text{ and } \norm{\bs{d}(v_\mathrm{true})-\bs{d}(v_\mathrm{best}^r)}/\norm{\bs{d}(v_\mathrm{true})} \simeq 0.1\%.
\end{equation} 
In other words, the component of $\log v_\mathrm{true}$ that is orthogonal to the span of the reduced basis has only a weak impact on the traveltimes.

\subsection{Inference results}\label{sec:tomo-online}
A \gls{pc} order of $15$ have been selected for the surrogate model of $A^{\pm 1/2}\Sigma^{\pm 1/2}(l)$ and $r\log A+ \mathrm{logdet}_\Sigma(l)$ ensuring a \gls{rrmse} around $10^{-3}$. Regarding the traveltimes, a level $3$ sparse grid with $15,135$ points (dimension $r+1=21$), requiring $75,675$ solves of the eikonal equation~\eqref{eq:eikonal} because of the $5$ sources, has been used with a \gls{rrmse} of $1\%$. 

\begin{figure}[!ht]
  \centering
    \centering
    \subfloat{
    \includegraphics[height=3.6cm]{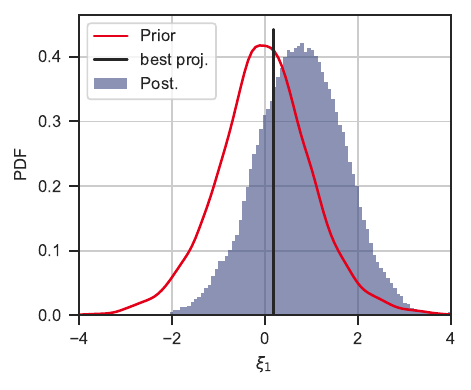}
    }\hspace{0.05cm}
    \subfloat{
      \includegraphics[height=3.6cm]{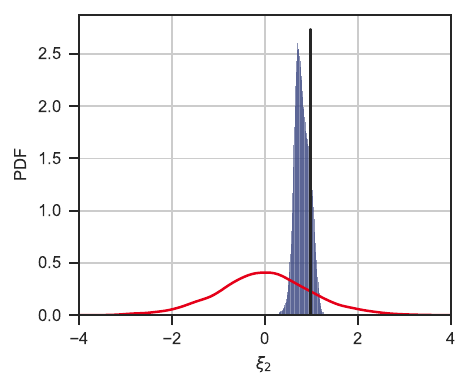}
    }
    \subfloat{ 
    \includegraphics[height=3.6cm]{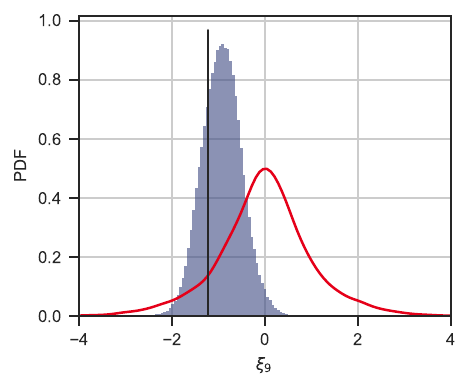}
    }\\\vspace{0.05cm}
    \subfloat{ 
    \includegraphics[height=3.6cm]{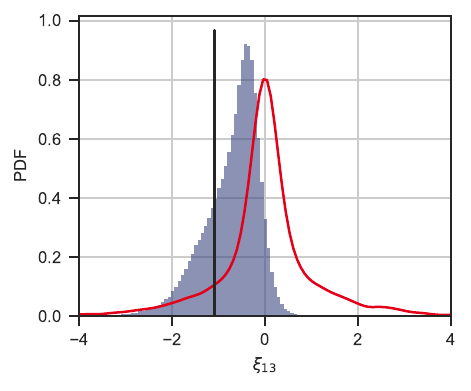}
    }\hspace{0.05cm}
    \subfloat{ 
      \includegraphics[height=3.6cm]{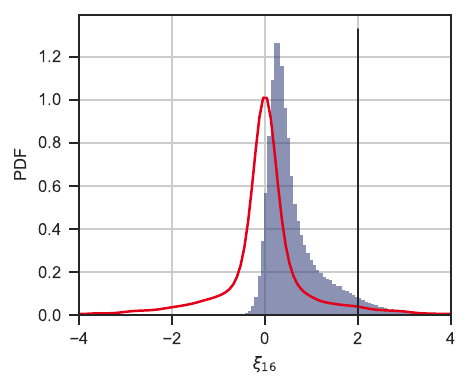}
    }\hspace{0.05cm}
    \subfloat{
      \includegraphics[height=3.6cm]{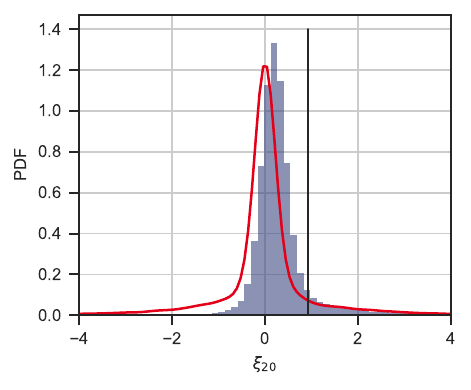}
    }\\\vspace{0.05cm}
    \subfloat{
    \includegraphics[height=3.6cm]{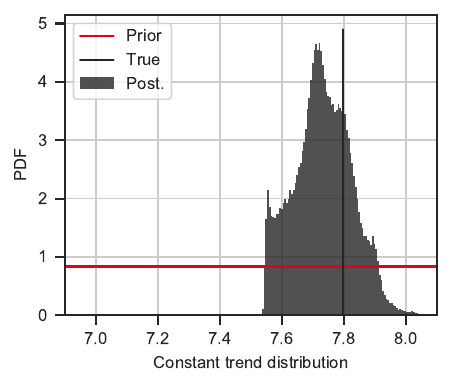}
    }\hspace{0.05cm}
    \subfloat{
      \includegraphics[height=3.6cm]{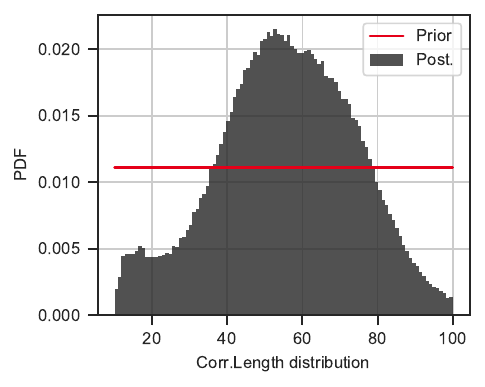}
    }\hspace{0.05cm}
    \subfloat{
      \includegraphics[height=3.6cm]{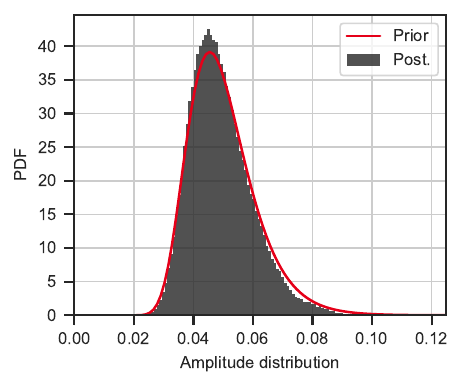}
    }
    \caption[Marginal posterior distributions of the coordinates and parameters.]{\acrshort{st} case - Marginal posterior distributions of some coordinates $\bxi$ and (hyper)parameters. Best projection coordinates are $g_\mathrm{kl}$ coordinates (sparse noised observations not taken into account).}\label{fig:mcmc-params}
\end{figure}

With the surrogate, a \gls{mcmc} chain burn-in of $10^6$ steps is performed, where the $(r+4)^2$ \gls{mcmc} covariance proposal is adapted every $5{\times}10^4$ steps. The acceptance rate is converging around $24\%$ at the end of the burn-in. After the burn-in, $5{\times}10^6$ steps are performed. The analysis shows an effective sample size around $6,000$ using the method proposed in~\cite{vats2019}. Figure~\ref{fig:mcmc-params} depicts the posterior marginal histograms of the coordinates together with their respective prior. The first coordinates, except for $\xi_1$, are rather well inferred with narrow distributions containing the projection coordinates while the last coordinates stay near their prior. The posterior of the constant trend concentrates.
 The first coordinate $\xi_1$ and the constant trend $c$ are difficult to distinguish since $\ol{u}_1$ is nearly constant. This explains why the first coordinate is not peaked in contrast to $\bxi_{2-10}$. The posterior of the hyperparameters are also plotted on Fig.~\ref{fig:mcmc-params}. We remark that the observations are weakly informative for the prior variance whose posterior distribution remains close to its prior.
The posterior distribution of the correlation length is more concentrated in the middle of the prior range and is asymetrical with a higher probability for low values. The fact that the posterior distribution of $l$ is not peaked highlights the interest of exploring the hyperparameters space to propose a large variety of plausible fields. Furthermore, we computed the probability of the true log field to be a Gaussian process realization with autocovariance function $k(\bq)$. The best hyperparameters couple is $(0.04,28)$ which is broadly consistent with the posterior distributions in spite of the noise level, that can affect these best values~\cite{musolas2021}. The noise level $\sigma_{\varepsilon}$ is correctly inferred with a distribution peaking at $1.8\mathrm{ms}$ with a standard deviation of $0.1\mathrm{ms}$, consistent with the true value of $2\mathrm{ms}$.
\begin{figure}[!ht]
  \centering
  \subfloat{
    \includegraphics[width=0.44\textwidth]{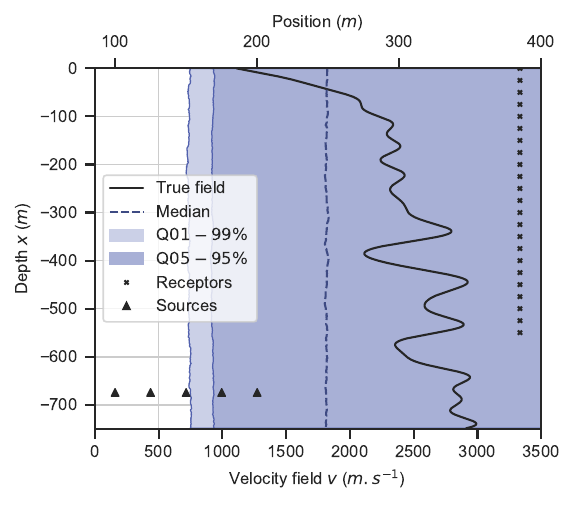}
  }\hfill
  \subfloat{
   \includegraphics[width=0.44\textwidth]{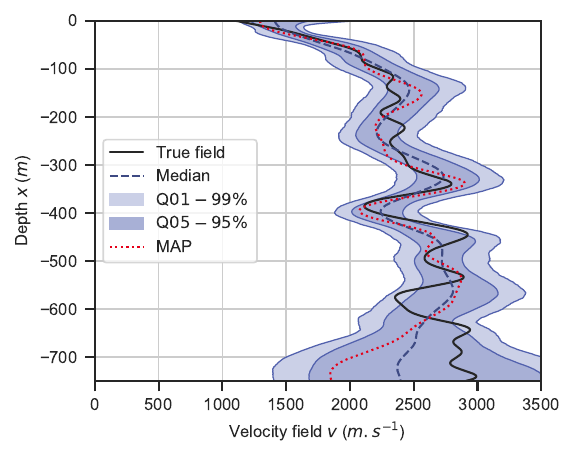}
  }
  \caption[Velocity field distribution.]{\acrshort{st} case - Distribution of the velocity field $f$ prior (left) and posterior (right). }\label{fig:mcmc-field}
\end{figure}

We now present our procedure to obtain the posterior distribution of the velocity field. To take into account the uncertainty due to the truncation, we add the $K$ next terms in the approximation as done in~\cite{meles2022}. \textit{A priori},
\begin{equation}
  g(\bx) = c + \suml_{i = 1}^r \ol{\lambda}_i^{1/2}\ol{u}_i(\bx) \xi_i + \suml_{i = r+1}^{r+K} \ol{\lambda}_i^{1/2}\ol{u}_i(\bx) X_i,
\end{equation}
with the $R\coloneqq r+K$ coordinates $(\bxi, \bs{X}) \sim \mc{N}(0,\Sigma(\bq))$, where $\Sigma(\bq) \in \R^{R\times R}$ is defined as in Eq.~\eqref{eq:sigma-def} and is partitioned as follows 
\begin{equation}
  \Sigma(\bq) = \begin{pmatrix} \Sigma_{rr}(\bq) & \Sigma_{rK}(\bq) \\ \Sigma_{Kr}(\bq) & \Sigma_{KK}(\bq) \end{pmatrix}.
\end{equation}
Thus, $\bs{X}|\bxi \sim \mc{N}(\ol{\mu}, \ol{\Sigma}(\bq))$, where 
\begin{equation}
  \ol{\mu} = \Sigma_{Kr}(\bq)\Sigma_{rr}(\bq)^{-1}\bxi \text{ and } \ol{\Sigma}(\bq) = \Sigma_{KK}(\bq) - \Sigma_{Kr}(\bq)\Sigma_{rr}(\bq)^{-1}\Sigma_{rK}(\bq).
\end{equation}
In practice, we use $R= r+K=20+60=80$. Since generating $\bs{X}$ knowing $\bxi$ and $\bq$ requires the construction of $\Sigma_{Kr}(\bq)$ and $\Sigma_{KK}(\bq)$ and the decomposition of $\ol{\Sigma}(\bq)$, the cost can be significant. Subsamples of the \gls{mcmc} chains are extracted and their velocity field are generated randomly. In practice, one sample over $100$ is selected.
The median, MAP and two quantiles range of the field posterior and prior distributions
are represented on Fig.~\ref{fig:mcmc-field}. The \gls{map} is the sample from the chain that maximizes the likelihood.
The median and the \gls{map} globally capture the main variations of the true field. 
The confidence intervals are quite tight (especially at the top of the domain) and remain consistent with the true field.
Larger uncertainties remains at the bottom of the domain because of the lack of observations there.

\subsection{Comparison for different fields}\label{sec:tomo:comp-field}
\begin{figure}[!ht]
  \centering 
  \subfloat[Posterior field distribution for $v_\mathrm{LW}$]{
    \includegraphics[width=0.44\textwidth]{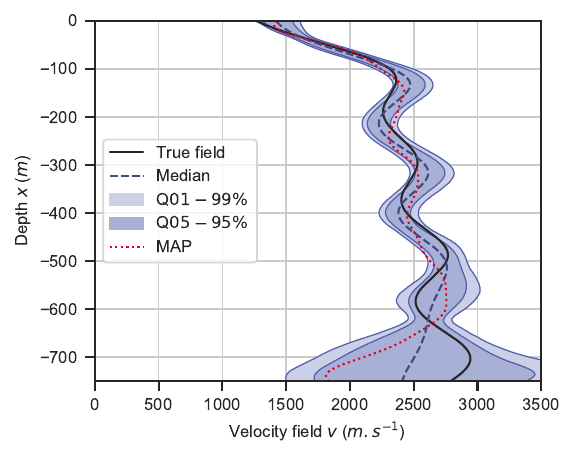}
  }\hfill 
  \subfloat[Posterior corr.length distribution for $v_\mathrm{LW}$]{
    \includegraphics[width=0.44\textwidth]{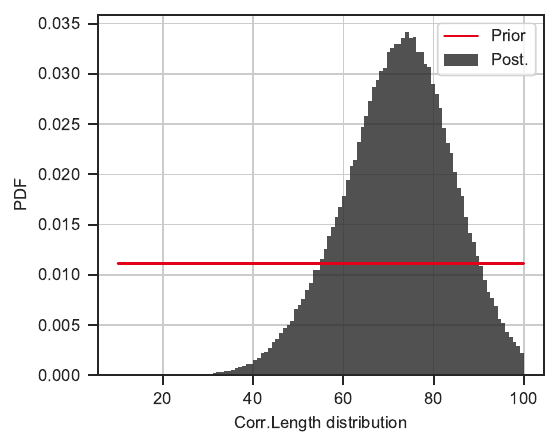}
  }\\
  \subfloat[Posterior field distribution for $v_\mathrm{SW}$]{
    \includegraphics[width=0.44\textwidth]{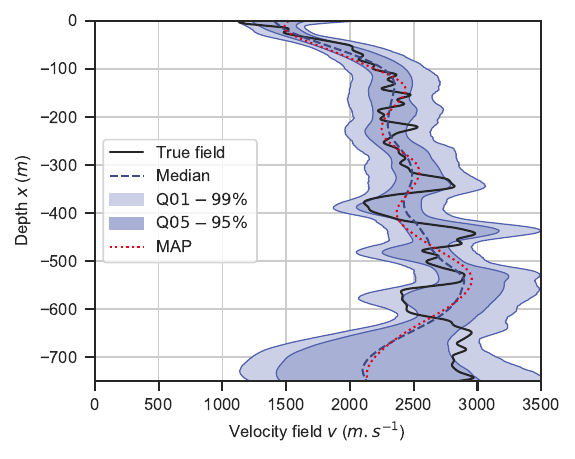}
  }\hfill 
  \subfloat[Posterior corr.length distribution for $v_\mathrm{SW}$]{
    \includegraphics[width=0.44\textwidth]{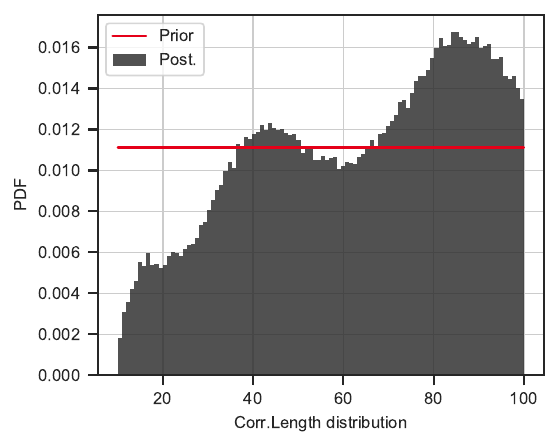}
  }
  \caption{\acrshort{st} case - Inference results for two fields with various wavelengths. Posterior distributions are obtained adding $K=60$ terms after the inference.}
  \label{fig:inference-diff-fields}
\end{figure}

To evaluate the robustness of the \gls{com} method, we repeat the inference for two versions of $v_\mathrm{true}$: 
\begin{inparaenum}[i)]
 \item $v_\mathrm{LW}$ obtained by smoothing $v_\mathrm{true}$ to remove short wavelengths, and 
 \item $v_\mathrm{SW}$ obtained by perturbing $v_\mathrm{true}$ with short wavelengths.
 \end{inparaenum} 
Posterior fields and correlation lengths distributions are shown on Fig.~\ref{fig:inference-diff-fields}.
The posterior uncertainties depend on the wavelengths of the field: they are rather smooth and tight for $v_\mathrm{LW}$ and less regular and larger for $v_\mathrm{SW}$. This behaviour is consistent with the modified true field shapes. 
Further, the correlation length distributions are impacted in a coherent way: the small correlation lengths are much less likely for $v_\mathrm{LW}$ than for $v_\mathrm{SW}$ whereas for $v_\mathrm{SW}$ the whole prior range remains plausible. 
Clearly, $v_\mathrm{SW}$ is more difficult to infer than $v_\mathrm{LW}$ because of its broader range of samples profile.
The inference of short wavelengths structures would require a denser arrangement of sources and receivers and possibly more information observations with incidence angle at arrival and/or full waveform analysis. 

\subsection{Comparison with fixed hyperparameters inference}\label{sec:comp-param}
To assess the interest of the \gls{com} method, we point out the limitations of using bases with fixed correlation lengths to infer $v_\mathrm{LW}$ and $v_\mathrm{SW}$. Two bases with different correlation lengths are built, 
\begin{inparaenum}[i)]
  \item $l=10$ to capture short wavelengths, and 
  \item $l=80$ to focus on long wavelengths.
\end{inparaenum}
For the two bases, an inference is led considering only the standard deviation noise $\sigma_\varepsilon$, the constant trend $c$ and the amplitude $A$ as hyperparameters, and inferring the coordinates in the basis parametrized with $l=10$ or $l=80$, with a standard Gaussian prior.

\begin{figure}[!ht]
  \centering 
  \subfloat[$v_\mathrm{LW}$, $l=10$]{
    \includegraphics[width=0.42\textwidth]{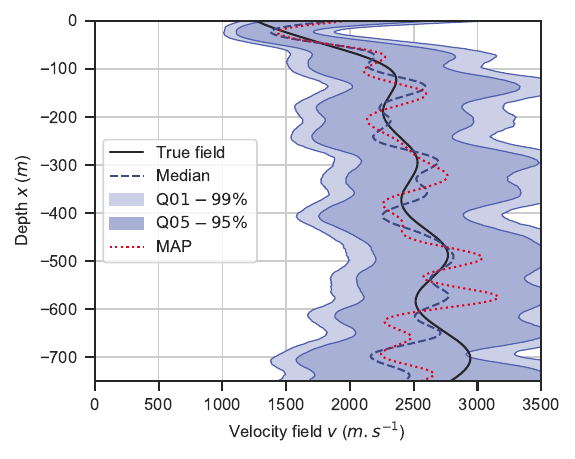}
  }\hfill 
  \subfloat[$v_\mathrm{LW}$, $l=80$]{
    \includegraphics[width=0.42\textwidth]{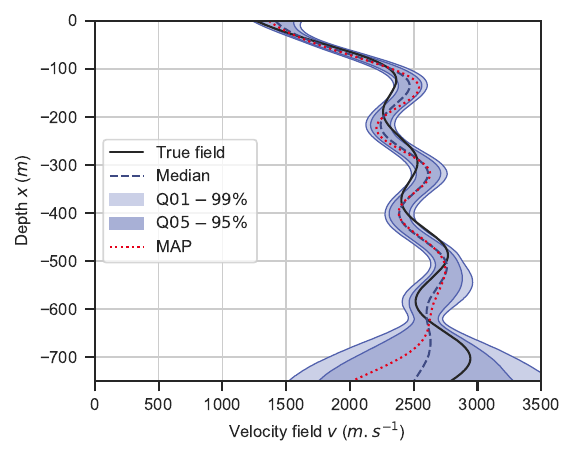}
  }\\
   \subfloat[$v_\mathrm{SW}$, $l=10$]{
    \includegraphics[width=0.42\textwidth]{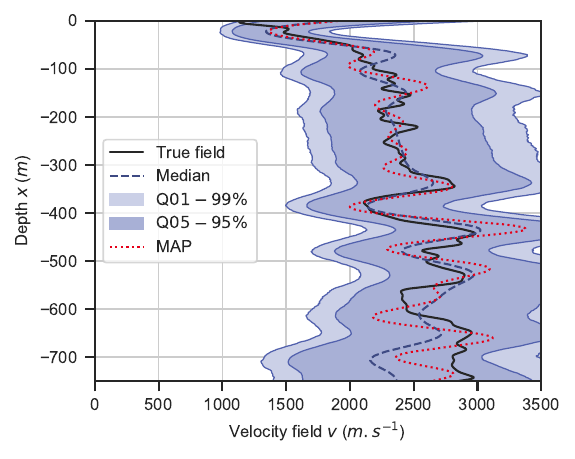}
  }\hfill 
  \subfloat[$v_\mathrm{SW}$, $l=80$]{
    \includegraphics[width=0.42\textwidth]{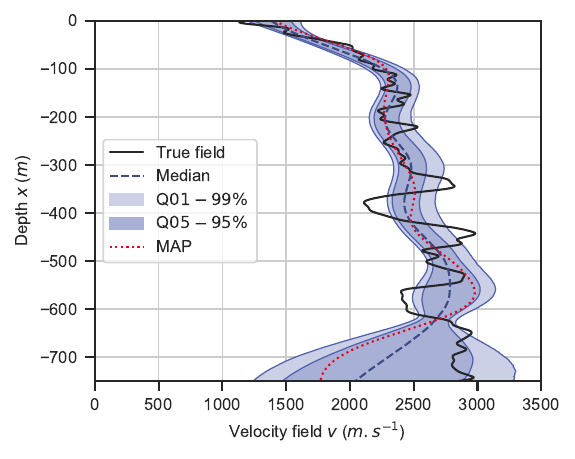}
  }
  \caption{\acrshort{st} case - Posterior field distribution for different field and inference bases. $r=20$, $K=60$ coordinates are added after the inference.}\label{fig:inference-field}
\end{figure}

Figure~\ref{fig:inference-field} displays the posterior field confidence interval obtained by summing $r=20$ (and $K=60$) coordinates. We see that a larger correlation length reduces the posterior uncertainties where the true field is well correlated. However, large $l$ values are unable to capture the small structures of $v_\mathrm{SW}$. On the contrary, smaller correlation lengths can account for small structures of $v_\mathrm{SW}$ but without reducing much the prior uncertainties. 
This result can be explained by the unsufficient information to learn the short wavelength features and the number of terms $r=20$ used. Comparing Fig.~\ref{fig:inference-diff-fields} and Fig.~\ref{fig:inference-field} shows the importance of allowing the inference of the covariance hyperparameters, since the hyperparameters choice is impacting the posterior distribution a lot. In turn, the exploration of the hyperparameters space encountered in the \gls{com} method improves the estimation of the field.

\section{Application to a 2D groundwater flow problem}\label{sec:groundwater}
The \gls{com} method is finally applied to a two-dimensional case with more hyperparameters to illustrate the approach's feasibility on a more realistic problem. The proposed test case introduced in Section~\ref{subsec:groundwater:case} is a source-sink groundwater flow problem inspired from~\cite{cui2014}. The inference results discussed in Section~\ref{subsec:groundwater:res} show that the \gls{com} method is able to deal with a parametrization adapted for anisotropic permeability fields. 

\subsection{Case presentation}\label{subsec:groundwater:case}
\begin{figure}
    \centering 
    \includegraphics[height=0.2\textheight]{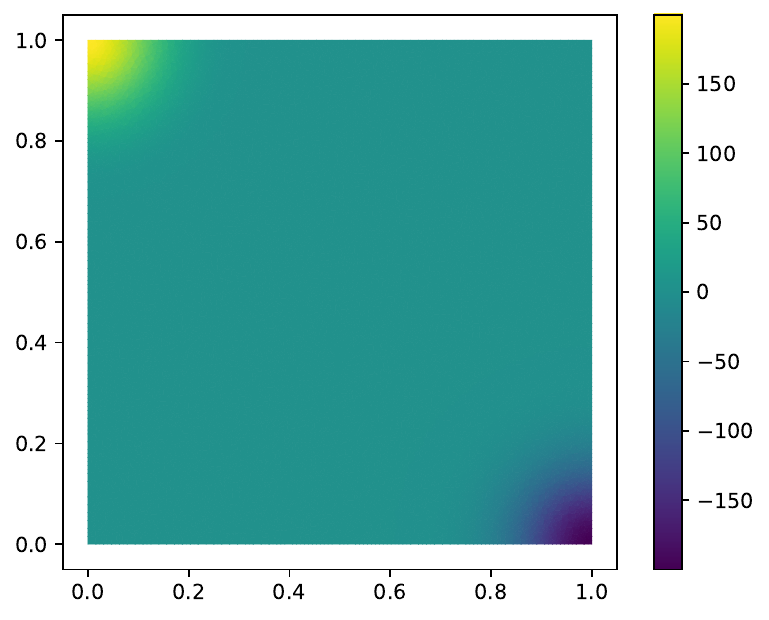}\hfill \includegraphics[height=0.2\textheight]{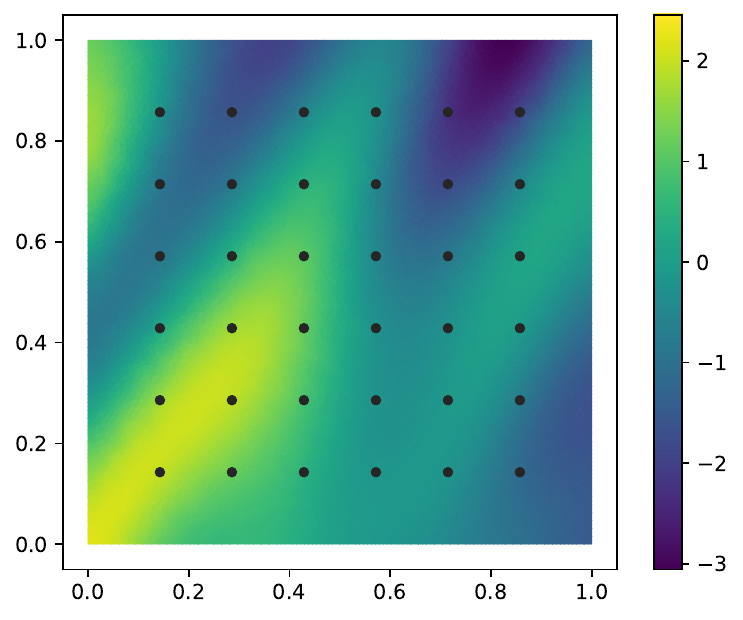} \hfill \includegraphics[height=0.2\textheight]{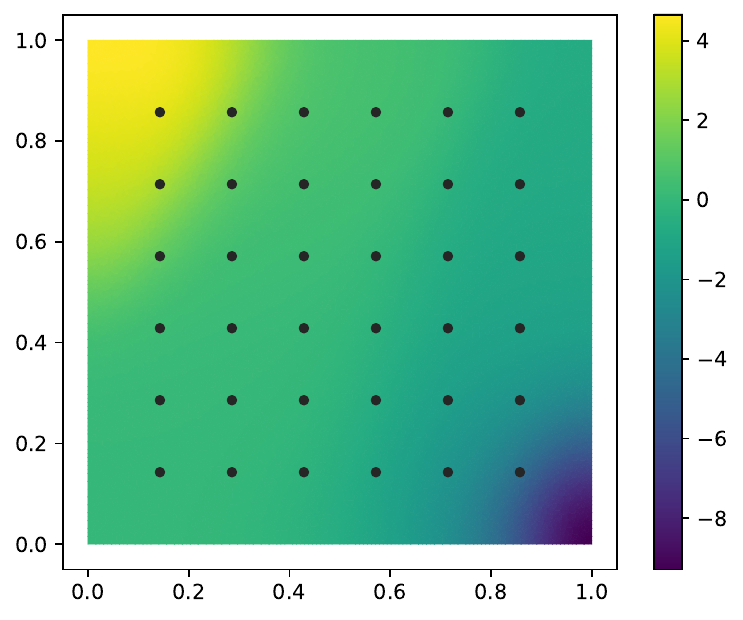} 
    \caption{Case presentation: (left) Source term $f$, (middle) True log-field $\log \kappa_\mathrm{true}$, (right) True solution $u_\mathrm{true}$. Black circles are observation locations.}\label{fig:groundwater:case-presentation}
\end{figure}

A steady-state groundwater flow problem with impermeable boundary condition is considered, 
\begin{equation}
    \left\{
\begin{alignedat}{3}\label{eq:groundwater:forward-pb}
    -\nabla \cdot (\kappa(\bx) \nabla u(\bx)) &= f(\bx), &&\quad \bx \in \Omega=[0,1]^2, \\
    \nabla u(\bx) \cdot \bs{n} &= 0, &&\quad \bx \in \partial \Omega \backslash \bs{0}, \\
    u(\bs{0}) &= 0, 
\end{alignedat}
    \right.
\end{equation}
where $u$ is the pressure field and $\kappa$ is the permeability field. We set $u(\bs{0}) = 0$ to ensure the well-posedness of the problem. The right-hand side term $f$, plotted in Fig.~\ref{fig:groundwater:case-presentation} (left), consists of a source in the top-left corner and a sink in the bottom-right corner, with rates following a squared-exponential shape centered on the corners and having a characteristic width of $0.1$. Given a realization of $\kappa$, the forward problem is solved by the $\mathbb{P}_2$ triangular finite element method using a mesh generated using~\cite{shewchuk1996}.

The objective is to infer the permeability field $\kappa$. The prior of $\log \kappa$ is a Gaussian distribution with an anisotropic squared exponential autocovariance function,
\begin{equation}\label{eq:groundwater:field-param}
\log \kappa \sim \mc{N}(0,k) \quad \text{with} \quad \left\{ \begin{aligned}
    &k(\bx,\by) = A\exp\left(\frac{-(\bx-\by)^\top K (\bx-\by)}{2}\right), \quad \forall \ \bx, \by \in \Omega, \\
    &K = R(\theta)\begin{pmatrix} 1/l_1^2 & 0\\ 0 & 1/l_2^2\end{pmatrix}R(\theta)^\top,
\end{aligned}\right.
\end{equation}
where $R(\theta)$ denotes the $2$D rotation matrix with orientation $\theta$ and $l_1$ and $l_2$ are the two correlation lengths characterizing the anisotropy. The four autocovariance hyperparameters priors are 
\begin{equation}
    A\sim \mathrm{IG}(3,14), \quad l_1 \sim \mc{U}(0.1,0.6), \quad l_2 \sim \mc{U}(0.1,0.6), \quad \theta \sim \mc{U}(0,\pi/2).
\end{equation}
The true field $\kappa_\mathrm{true}$ used to generate the observations is a particular realization of the random process described in Eq.~\eqref{eq:groundwater:field-param} with $\{A^\star,\ l_1^\star,\ l_2^\star,\ \theta^\star\} = \{1.3,\ 0.4,\ 0.15,\ \pi/3\}$. This realization is shown in the center plot of Fig.~\ref{fig:groundwater:case-presentation}, along with its associated pressure field $u_\mathrm{true}$ (right plot) computed with a fine uniform triangulation of the domain with $27,763$ elements. The observations correspond to the pressure field $u_\mathrm{true}$ measured over a uniform grid of $6\times 6$ sensors (see Fig.~\ref{fig:groundwater:case-presentation}).  An \gls{iid} Gaussian noise of standard deviation $\sigma_\varepsilon = 0.02$ (around $1\%$ of the average magnitude of the true values) is added to play the role of measurement errors. 

\begin{figure}
    \includegraphics[width=0.32\textwidth]{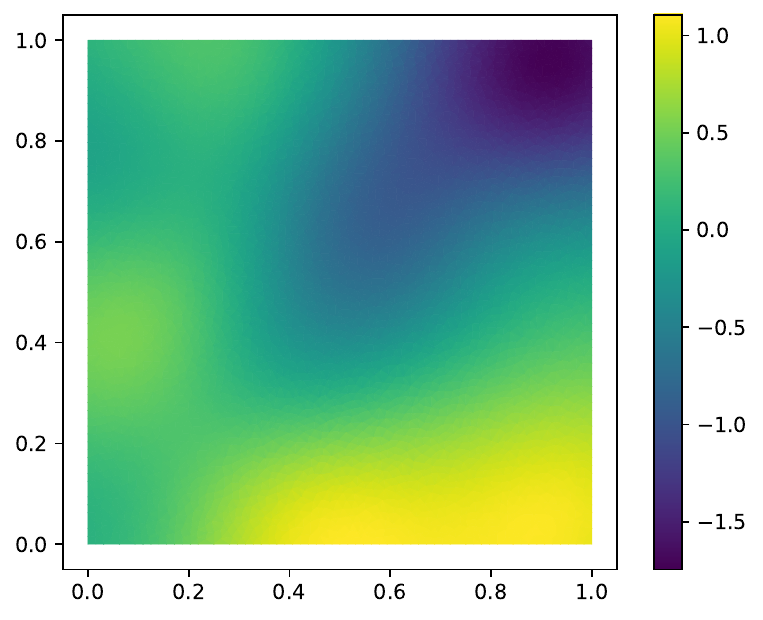}\hfill \includegraphics[width=0.32\textwidth]{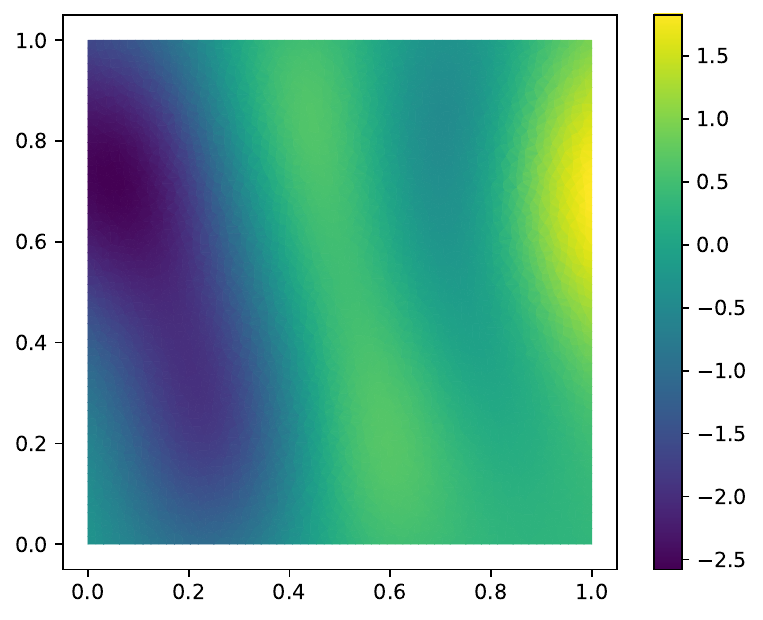} \hfill \includegraphics[width=0.32\textwidth]{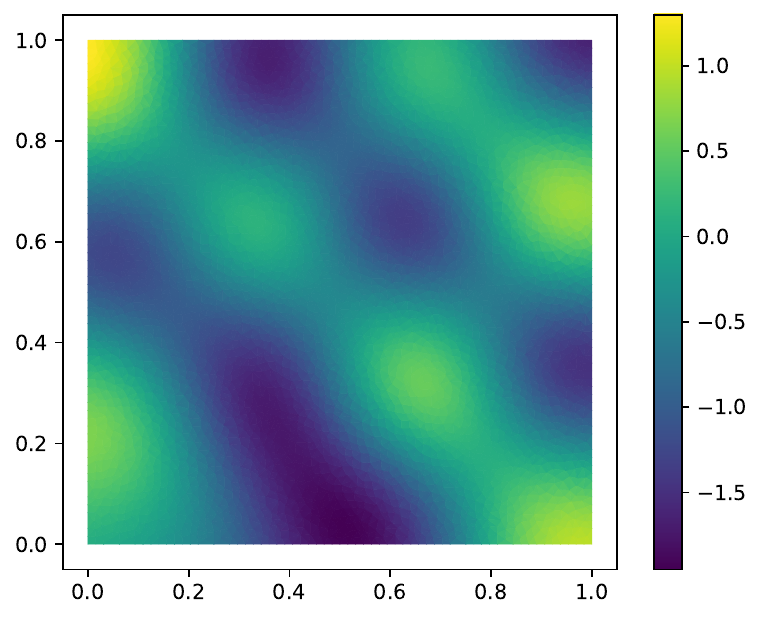}
    \caption{Few samples $\log \kappa$.}\label{fig:groundwater:prior}
\end{figure}

The truncated \gls{kl} decomposition used for the inference has $r=20$ modes which capture $97\%$ of the prior variance. The inference uses a coarser mesh composed of $7,875$ elements to introduce a model error term. Using this parametrization, some samples of $\log \kappa$ drawn from its prior are shown in Fig.~\ref{fig:groundwater:prior} to highlight the rich structure of the~\textit{a priori}.

\subsection{Results}\label{subsec:groundwater:res}
The parameters to infer are $r=20$ reduced coordinates, the amplitude $A$, the two correlation lengths $l_1$ and $l_2$, the orientation $\theta$ as well as the noise standard deviation $\sigma_\varepsilon$, leading to a sampling space of dimension $25$.
As for the previous test cases, \gls{pc} surrogates are built for $\Sigma(\bq)^{\pm 1/2}/A^{\pm 1/2}$ and $\mathrm{logdet}_\Sigma(\bq) - r\log A$ where $\bq = \{l_1,\ l_2,\ \theta\}$ (see Section~\ref{subsec:acceleration:surrogate-prior}). A \gls{rrmse} lower than $1\%$ is achieved with a $12$th order \gls{pc} expansion (total degree truncation) fitted using a Latin hypercube sampling of size $1,000$. The $\mathbb{P}_2$ finite element method used to solve Pb.~\eqref{eq:groundwater:forward-pb} is fast enough to not require the use of a forward model surrogate.

\begin{figure} 
    \includegraphics[width=0.24\textwidth]{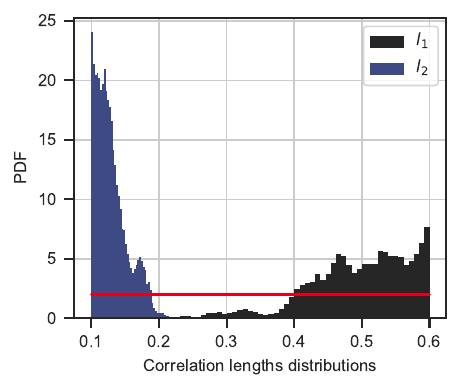}
    \includegraphics[width=0.24\textwidth]{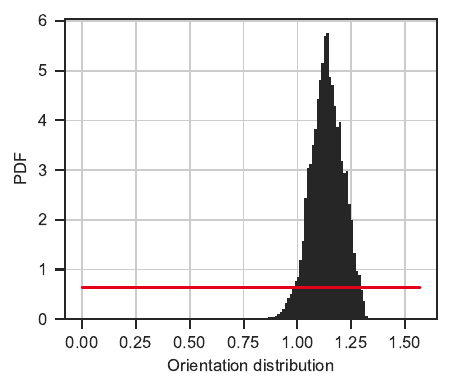}
    \includegraphics[width=0.24\textwidth]{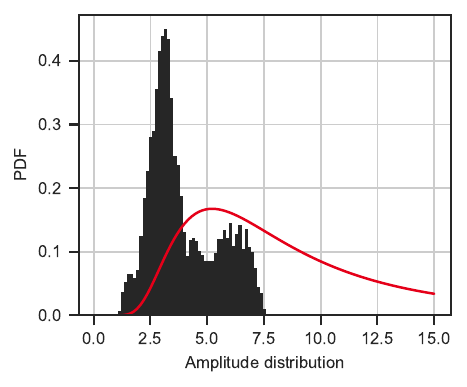} 
    \includegraphics[width=0.24\textwidth]{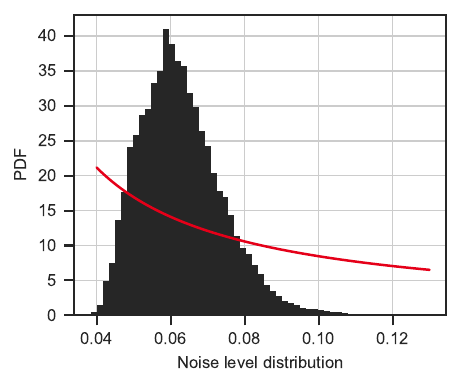}
\caption{Marginal distributions of kernel hyperparameters (left to right) first and second correlation lengths $l_1$, $l_2$, orientation $\theta$, amplitude $A$, noise level $\sigma_\varepsilon$}\label{fig:groundwater:hyperparameters}
\end{figure}

The posterior marginals of the five hyperparameters are reported in Fig.~\ref{fig:groundwater:hyperparameters}. The marginal of the orientation peaks around $\pi/3$ while the marginal of the two correlation lengths are distinct. The extent of the marginal of $A$ is sensibly reduced compared to its prior distribution. These marginals are somehow consistent with the particular values of hyperparameters used for generating the true field. Finally, the marginal of the noise level shows likely values slightly higher than the value used for generating the observations, as expected due to a non-vanishing model error induced by using a coarser mesh and a truncation of the field representation for the inference.

\begin{figure} 
    \includegraphics[width=0.24\textwidth]{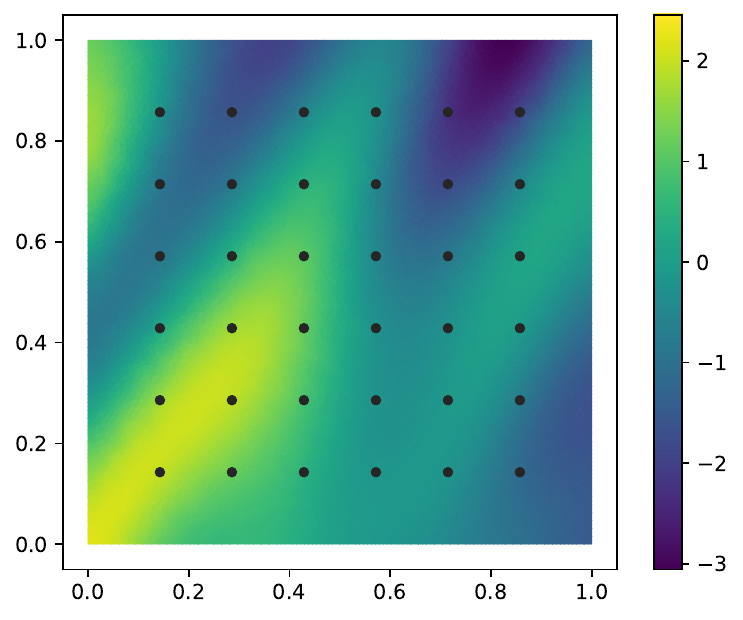}
    \includegraphics[width=0.24\textwidth]{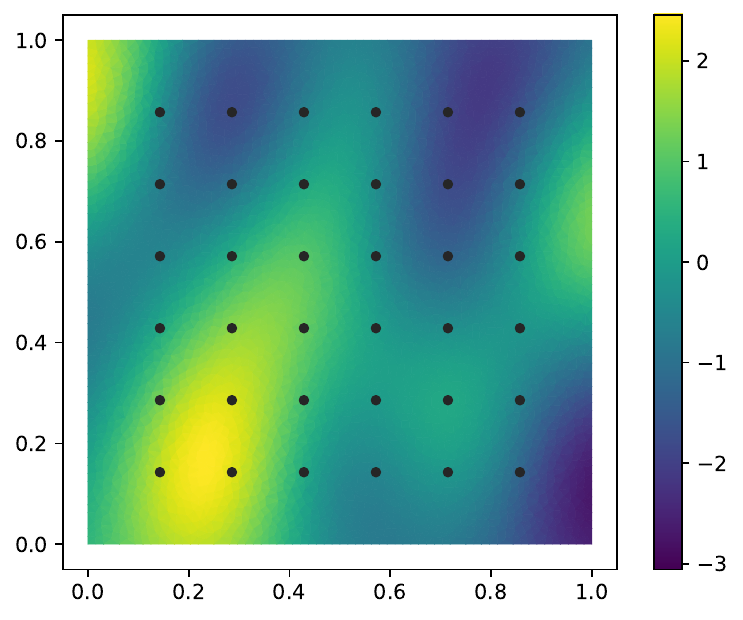}
    \includegraphics[width=0.24\textwidth]{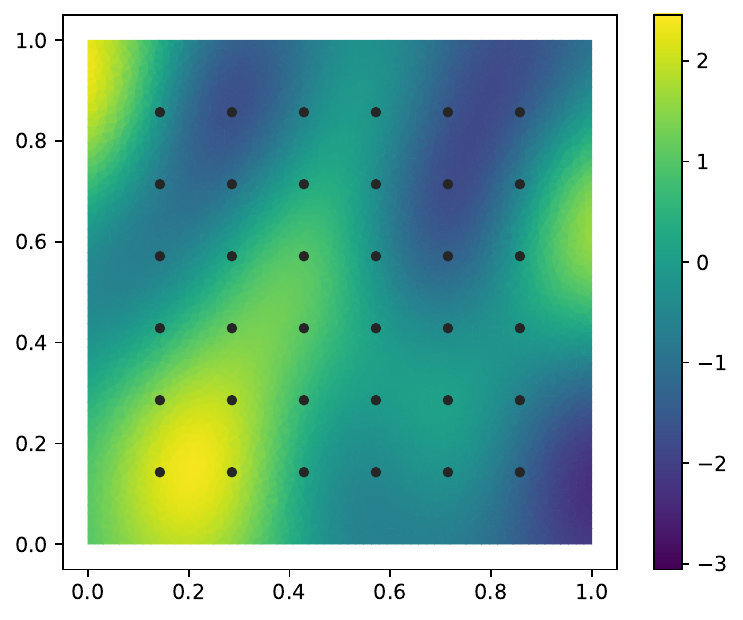}
    \includegraphics[width=0.24\textwidth]{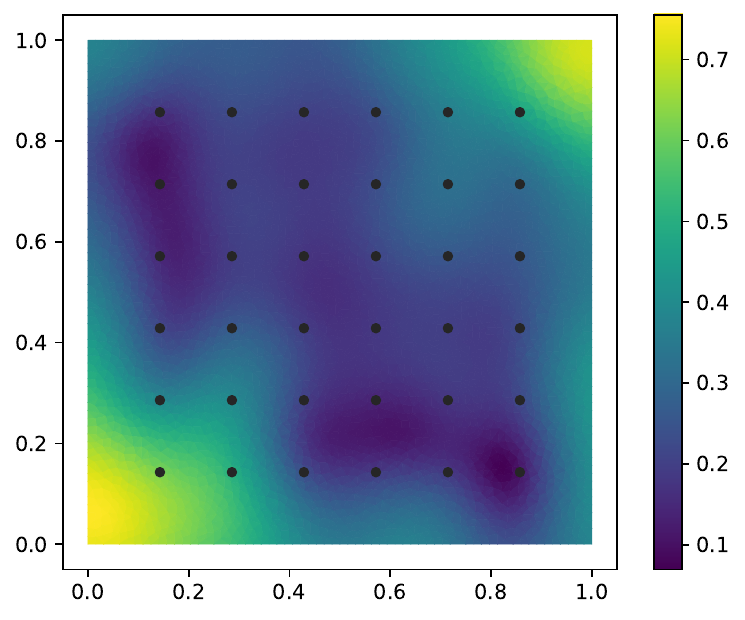}
\caption{(second plot) \gls{map} (third plot) mean and (fourth plot) standard deviation of the posterior $\log \kappa$ distribution. The true log-field is recalled (first plot) in order to facilitate the comparison.}\label{fig:groundwater:fielddistrib}
\end{figure}

Figure~\ref{fig:groundwater:fielddistrib} shows the true log field (first plot), its \gls{map} estimate (second plot), posterior mean (third plot), and standard deviation (fourth plot). It is seen that the \gls{map} and the posterior mean field are close to the true field. Further, the posterior uncertainty, as measured by the standard deviation, is higher in the corners far from the observations and the injection and pumping areas.

\begin{figure}
\centering 
\includegraphics[width=0.4\textwidth]{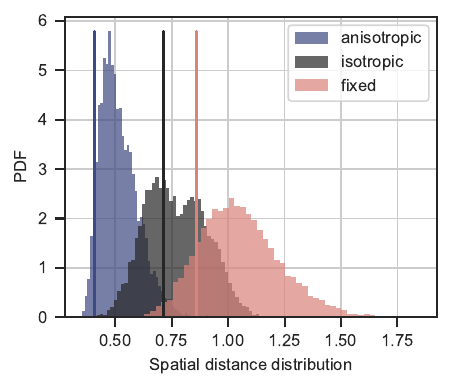}
\caption{Posterior distributions of the spatial distance between the true log-field and posterior realizations $\left|\left| \log \kappa_{true} - \log \kappa^{(i)} \right|\right|_{L^2(\Omega)}$ for the (blue) anisotropic (gray) isotropic (pink) fixed hyperparameters case. Distances between the true log-field and the \gls{map} log-fields are indicated with lines.}\label{fig:groundwater:comp-dist}
\end{figure}

To highlight the benefit of using a parametrized prior, we compared the inferred fields with two other approaches based on an isotropic covariance ($l_1=l_2$, no orientation). The first relies on the \gls{com} method with parameters $\{A,l\}$ and the second uses a \gls{kl} parametrization with a fixed correlation length $l=0.3$ and amplitude $A=1.5$. In all cases, we use $r=20$ reduced coordinates in the inference.
We assess the inference quality by reporting the posterior distribution of the norm of the difference between the inferred log-fields and the true log-field. These distribbutions, reported in Fig.~\ref{fig:groundwater:comp-dist}, show that a richer parametrization of the prior significantly reduces the posterior distance to the true field. The error on the \gls{map} estimates, shown with the vertical lines, exhibits consistent behavior. These results demonstrate the importance of effectively treating hyperparameters to enable more general priors.

\section{Conclusion}\label{sec:ccl}
In this paper, we developed and implemented a novel approach to efficiently deal with hyperparameters in Bayesian inference of Gaussian fields. Specifically, the approach relies on a representation of the field in a fixed reduced basis, carefully selected, and the prior's covariance of the reduced coordinates is seen as a function of the hyperparameters.
The cornerstone of our approach is a change of measure designed to avoid the dependency of the likelihood on the hyperparameters. From a methodological viewpoint, the change of measure 
yields a posterior that depends smoothly on the hyperparameters and reduced coordinates, compared to alternative approaches such as the change of coordinates which can exhibit discontinuities.

This smoothness enables to use surrogate models to decrease the computational cost of the MCMC sampling.
In practice, the \gls{com} results in a hierarchical Bayes formulation with a closed-form expression for the hyperparameters-dependent prior: it allows the efficient sampling of the joint distribution with an auxiliary variable. The method was demonstrated on three test cases with various complexities and geometries. 
The numerical experiments demonstrate that the hyperparameters space exploration is 
\begin{inparaenum}[i)]
    \item highly valuable because it provides a better estimation of the field uncertainties and
    \item computationally tractable thanks to the \gls{com} formulation combined with surrogate models.
\end{inparaenum}

Further works should focus on the extension of \gls{com} methods to non-Gaussian random processes with possibly non stationnary prior autocovariance functions. In that case, the change of measure formulation needs to be generalized to account for possible non Gaussian dependencies in the priors. For seismic tomography applications, it would be interesting to add other types of observations, such as the direction of propagation, the maximum amplitude or even the full waveform, to improve the estimation of the reduced coordinates with higher indices. More generally, the extension to three-dimensional cases is not conceptually limited by the \gls{com} method but by the number of \gls{kl} modes required to accurately represent the field of interest. Current work focuses on the development of adaptive dimension reduction techniques to tackle this issue.

\bibliographystyle{elsarticle-num}
\bibliography{COMarticle}

\appendix 
\section{Marginal distribution of the coordinates}\label{appendix:xi-law}
In Section~\ref{sec:acceleration}, the approximation $\bxi\sim\mc{N}(0,\mathrm{I}_r)$ is used for the \gls{pc} surrogate construction. This choice is justified in the following. The conditional distribution of the coordinates $\bxi$ according to the hyperparameters $\bq$ writes 
\begin{equation}
  \bxi_{|\bq} \sim \mc{N}(0,\Sigma(\bq))\coloneqq \mc{N}_0(\bq), \quad \bq \sim \pi_\mathbb{H}
\end{equation}
The marginal distribution of the coordinates $\bxi$ over the hyperparametric domain is a compound probability distribution. Note that this distribution has no reason to be Gaussian and its exact shape depends on the distribution of the hyperparameters.
Here, we want to find the normal distribution $\mc{N}_1(\mu, C)$ such that it approximates the best $\bxi_{|\bq}$. That is to say, we want $\mc{N}_1$ to be the closest normal distribution in average along $\bq$ for a given distance. The Kullback--Leibler divergence between two multivariate Gaussian distribution is
\begin{equation}
2D_{\mathrm{KL}}(\mc{N}_0(\bq)||\mc{N}_1) = \mu^\top C^{-1} \mu + \mathrm{Tr}(C^{-1}\Sigma(\bq)) - \log \frac{\mathrm{det}(\Sigma(\bq))}{\mathrm{det}(C)} - r,
\end{equation}
where $r$ is the dimension. Denoting $\ol{\Sigma}$ the expectation of $\Sigma(\bq)$ with respect to $\bq$, we have
\begin{equation}\label{eq:sig-mean}
  \ol{\Sigma}_{ij} = \Int_\mathbb{H} \Sigma(\bq)_{ij}\pi_\mathbb{H}(\bq) d\bq \stackrel{\eqref{eq:sigma-def}}{=} (\ol{\lambda}_i\ol{\lambda}_j)^{-1/2}\scal{\scal{\Int_\mathbb{H}k(\bq)\pi_\mathbb{H}(\bq)d\bq}{\ol{u}_i}}{\ol{u}_j} \stackrel{\eqref{eq:kref-def}}{=} (\ol{\lambda}_i\ol{\lambda}_j)^{-1/2}\scal{\scal{\ol{k}}{\ol{u}_i}}{\ol{u}_j} \stackrel{\eqref{eq:Bref-def}}{=} \delta_{ij}.
\end{equation}
By linearity and using Eq.~\eqref{eq:sig-mean}, we get 
\begin{equation}
  \E_\mathbb{H}(\mathrm{Tr}(C^{-1}\Sigma(\bq))) = \mathrm{Tr}(C^{-1}\E_\mathbb{H}(\Sigma(\bq))) = \mathrm{Tr}(C^{-1}\ol{\Sigma}) = \mathrm{Tr}(C^{-1}).
\end{equation}
So, the $\bq$-averaged Kullback--Leibler divergence simplifies as 
\begin{equation}
  \E_\mathbb{H} (2D_{\mathrm{KL}}(\mc{N}_0(\bq)||\mc{N}_1) ) = \mu^\top C^{-1} \mu + \mathrm{Tr}(C^{-1}) + \log \mathrm{det}(C) -\E_\mathbb{H} (\log \mathrm{det}(\Sigma(\bq))) - r.
\end{equation}
Since $C$ is symmetric positive definite, $\mu$ is set to zero to minimize the divergence. The covariance $C^\ast$ minimizes therefore 
\begin{equation}
  C^\ast = \argmin_{C\in \mc{S}^r_{+}} \left[ \mathrm{Tr}(C^{-1}) + \log \mathrm{det}(C)\right],
\end{equation}
where $\mc{S}^r_{+}$ denotes the set of the symmetric positive definite matrices of size $r\times r$. Furthermore, $C \in \mc{S}^r_{+}$ can be decomposed 
\begin{equation}
  C = QLQ^{-1},
\end{equation}
with $Q$ an unitary matrix and $L$ a positive diagonal matrix. Then,
\begin{equation}
  \log \mathrm{det}(C) = \suml_{i=1}^r \log L_{ii} \quad \text{and} \quad \mathrm{Tr}(C^{-1}) = \mathrm{Tr}(QL^{-1}Q^{-1}) = \suml_{i=1}^r L_{ii}^{-1}
\end{equation}
The minimization consists in finding the positive diagonal elements of $L$, $L^\ast_{11},\dots,L^\ast_{rr}$ such that 
\begin{equation}
  L^\ast_{11},\dots,L^\ast_{rr} = \argmin_{L_{11},\dots,L_{rr}} \left[ \suml_{i=1}^r L_{ii}^{-1}+\log L_{ii}\right].
\end{equation}
It is a sum of positive terms whose minimum is achieved when all the terms are minimal. It corresponds to $L_{ii} = 1$ for $i=1,\dots,r$, leading to 
\begin{equation}
  C^\ast = QL^\ast Q^{-1} = Q\mathrm{I}_r Q^{-1} = \mathrm{I}_r
\end{equation}
To conclude, the normal distribution that minimizes the $\bq$-averaged Kullback--Leibler divergence to $\mc{N}_0(\bq)$ is $\mc{N}(0,\mathrm{I}_r)$.

\end{document}